\documentclass[10pt,reqno]{amsart}

\usepackage{amsmath,amsthm,amssymb,tikz-cd,comment,hyperref,cleveref,nicematrix,enumitem,breqn,float,caption}

\allowdisplaybreaks
\theoremstyle{plain}
\newtheorem{thm}{Theorem}[section]
\newtheorem{lem}[thm]{Lemma}
\newtheorem{prop}[thm]{Proposition}

\newtheorem{thmalphabetintro}{Theorem}

\newtheorem{thmalphabetmaintext}{Theorem}

\theoremstyle{definition}
\newtheorem{defn}[thm]{Definition}
\newtheorem{exmp}[thm]{Example}

\newtheorem{rem}[thm]{Remark}
\newtheorem{conjecture}[thm]{Conjecture}
\newtheorem{question}[thm]{Question}

\newtheorem{remark and notation}[thm]{Remark and Notation}

\def \p {\mathbb{P}}

\def \z {\mathbb{Z}}

\def \q {\mathbb{Q}}
\def \k {\mathbb{K}}

\def \o {\mathcal{O}}
\def \i {\mathcal{I}}
\def \f {\mathcal{F}}

\def \H {\textnormal{H}}
\def \sing {\textnormal{sing}}

\def \leng {\operatorname{length}}

\def \proj {\operatorname{Proj}}
\def \codim {\operatorname{codim}}

\def \tor {\operatorname{Tor}}
\def \reg {\operatorname{reg}}

\def \maxdeg {\operatorname{maxdeg}}
\def \mon {\operatorname{mon}}

\def \conv {\operatorname{conv}}

\def \coloredone {\color{gray}{1}}
\def \coloredtwo {\color{gray}{2}}
\def \boldone {1}
\def \matinds[#1] {\color{gray!60}{#1}}
\def \entryzero {}

\def \rectangleEntryColor[#1][#2] {\rectanglecolor[HTML]{EEEEEE}{#1}{#2}}

\title[Surface counterexamples to Eisenbud-Goto conjecture]{Surface counterexamples to the Eisenbud-Goto conjecture}
\author{Jong In Han}
\address{Jong In Han, Department of Mathematical Sciences, Korea Advanced Institute of Science and Technology (KAIST), 291, Daehak-ro, Yuseong-gu, Daejeon, Republic of Korea}
\email{jihan09@kaist.ac.kr}

\author{Sijong Kwak}
\address{Sijong Kwak, Department of Mathematical Sciences, Korea Advanced Institute of Science and Technology (KAIST), 291, Daehak-ro, Yuseong-gu, Daejeon, Republic of Korea}
\email{sjkwak@kaist.ac.kr}

\thanks{J. I. Han was partially supported by Basic Science Research Program through the National Research Foundation of Korea(NRF) funded by the Ministry of Education(2019R1A6A1A10073887). S. Kwak was partially supported by the National Research Foundation of Korea(NRF) grant funded by the Korea government(MSIT) (No. 2021R1A2C1013851).}

\begin{document}
	\begin{abstract}
		It is well known that the Eisenbud-Goto regularity conjecture is true for arithmetically Cohen-Macaulay varieties, projective curves, smooth surfaces, smooth threefolds in $\mathbb{P}^5$, and toric varieties of codimension two.
		After J. McCullough and I. Peeva constructed counterexamples in 2018, it has been an interesting question to find the categories such that the Eisenbud-Goto conjecture holds.
		So far, surface counterexamples have not been found while counterexamples of any dimension greater or equal to 3 are known.
		
		In this paper, we construct counterexamples to the Eisenbud-Goto conjecture for projective surfaces in $\mathbb{P}^4$ and investigate projective invariants, cohomological properties, and geometric properties. The counterexamples are constructed via binomial rational maps between projective spaces.
	\end{abstract}
	
	\maketitle
	
	\section{Introduction}
	For a coherent sheaf $\f$ on a projective space $\p^r$ over a field $\k$, the sheaf $\mathcal{F}$ is \textit{$m$-regular} if
	\[
	\H^i(\p^r,\f(m-i))=0\text{ for all }i\geq 1.
	\]
	The minimum integer $m$ such that $\f$ is $m$-regular is called the \textit{regularity of $\f$} and denoted as $\reg \f$.
	For a projective variety $X\subseteq\p^r$ with the ideal sheaf $\i_X$, the regularity $\reg X$ of $X$ is defined as the regularity of $\i_X$.
	Let $S=\k[x_0,\cdots,x_r]$ be the polynomial ring of $\p^r$.
	For a graded $S$-module $M$, the regularity of $M$ is defined as
	{\small
		\[
		\reg M:=\max\{j\in\z\mid\beta_{i,j}(M)=\dim \tor_i(M,\k)_{i+j}\text{ is nonzero for some }i\}.
		\]}
	
	\noindent
	If $M$ is saturated, i.e. $M\cong \bigoplus_{d\in \z}\H^0(\p^r,\widetilde{M}(d))$, then $\reg M=\reg \widetilde{M}$. In particular, we have $\reg I_X=\reg \i_X$ where $I_X$ is the homogeneous ideal of $X$.
	
	The regularity of projective varieties is considered as an important projective invariant measuring the complexity of projective varieties. For any arithmetically Cohen-Macaulay variety $X$, one can see there is a bound of the regularity:
	\[
	\reg X\leq \deg X-\codim X+1.
	\]
	In 1983, Gruson, Lazarsfeld, and Peskine proved that the same bound of regularity holds for every projective curve $X\subseteq\p^r$ (\cite{GLP}). In 1984, Eisenbud and Goto suggested the famous Eisenbud-Goto conjecture that claims the same inequality holds for every projective variety (\cite{EG}). In 1987, Lazarsfeld proved that the conjecture holds for smooth surfaces over the fields of characteristic zero using the method of general outer projections and vector bundle techniques on the regularity (\cite{L}).
	For the case of most smooth threefolds, Z. Ran studied using a local differential geometric method (\cite{R}).
	For other smooth threefolds, the second author proved there is a weak bound
	\begin{equation}\label{loose_ineq}
		\reg X\leq \deg X-\codim X+2
	\end{equation}
	using Lazarsfeld's construction and Mather's theorem (\cite{K2}).
	For smooth threefolds in $\p^5$, the Eisenbud-Goto conjecture holds with sharp cases as proved in (\cite{K3}).
	For smooth projective varieties of larger dimension, it is known that
	\[
		\reg X\leq \deg X-\codim X+1+\frac{(n-2)(n-1)}{2}
	\]
	where $n=\dim X\leq 14$ by L. Chiantini, N. Chiarli, and S. Greco (\cite{CCG}).
	
	For singular varieties, not many things are known although there are some positive results for normal varieties with mild singularities.
	In 2015, W. Niu proved the conjecture for normal surfaces with rational, Gorenstein elliptic, and log canonical singularities (\cite{N1}). In 2022, W. Niu and J. Park proved the weak inequality (\ref{loose_ineq}) for threefolds with mild singularities (\cite{NP2}).
	
	In 2018, J. McCullough and I. Peeva published a surprising result that there are counterexamples to the Eisenbud-Goto conjecture. They showed the regularity of projective varieties is not bounded by any polynomial function of the degree. Also, their example shows there is a nondegenerate singular threefold in $\p^5$ of degree 31 and regularity 38 (\cite[Example 4.7]{MP1}). It is regular in codimension one as proved in \cite{MMM}. Recently, J. Choe constructed some counterexamples by using the method called unprojection (\cite{C2}).
	
	In spite of the counterexamples, the Eisenbud-Goto conjecture is still believed to hold for projective varieties with nice properties.
	Also as J. McCullough and I. Peeva stated in \cite[12.6]{MP2}, the following are interesting cases where the conjecture may hold:
	\begin{itemize}
		\item smooth projective varieties,
		\item projectively normal varieties,
		\item projective toric varieties,
		\item projective singular surfaces.
	\end{itemize}
	
	However, in this paper, we construct counterexamples of projective singular surfaces of codimension two to the Eisenbud-Goto conjecture. As far as we know, these are the first surface counterexamples to the conjecture. The counterexamples are constructed via binomial rational maps between projective spaces.
	
	The Eisenbud-Goto regularity conjecture for $X$ is equivalent to $X$ satisfying both of the following condition:
	\begin{itemize}
		\item (normality conjecture) $X$ is $k$-normal for all $k\geq d-e$;
		\item (regularity conjecture for $\o_X$) $\reg\o_X\leq d-e$.
	\end{itemize}
	The problem identifying the number $k$ such that $X$ is $k$-normal is considered as an important problem and studied from the classical era, especially by Castelnuovo (\cite{C1}). He showed the smooth curves are $k$-normal for all $k\geq d-2$ which implies the normality conjecture for smooth curves in $\p^3$. Then L. Gruson, R. Lazarsfeld, and C. Peskine showed the normality conjecture for all projective curves (\cite{GLP}). For curves and surfaces, the Eisenbud-Goto conjecture is equivalent to the normality conjecture. The regularity conjecture for $\o_X$ is proved for smooth projective varieties by A. Noma, the second author, and J. Park (\cite{N2}, \cite{KP}). Hence the Eisenbud-Goto regularity conjecture is equivalent to the normality conjecture also for smooth projective varieties.
	
	To state counterexamples, we construct a family $\{X_m\}_{m\geq 6}$ of surfaces in $\p^4$.
	We highly suspect that all members of this family are counterexamples to the regularity conjecture. 
	For a homogeneous ideal $I$, we denote the largest degree of an element in a minimal generating set of $I$ as $\maxdeg I$, which is less or equal to $\reg I$.
	For a projective variety $X$ with homogeneous ideal $I_X$, we define $\maxdeg X:=\maxdeg I_{X}$.
	
	\begin{thmalphabetintro}\label{mainthm_general_m_intro}
		Let $m\geq 6$ be an integer and $X_m$ be the image of the binomial rational map \[\varphi_m:\p^2=\proj\k[y_0,y_1,y_2]\dashrightarrow\p^4=\proj\k[x_0,\cdots,x_4]\] given by the linear system \[(\o_{\p^2}(m),y_1^{m-1}y_2,
		y_0^{m-1}y_1,
		y_1^m+y_1^{m-3}y_2^3,
		y_0^m,
		y_2^{m-1}y_0).\]
		Then $X_m\subseteq\p^4$ is a surface of degree $m^2-m+3$.
		In addition, suppose $m=6k$ for some integer $k\geq 1$ and define the square matrices $L(k)$ and $W(k)$ as in \Cref{defn_N_L}. If $L(k)$ and $W(k)$ are both invertible, then $\maxdeg X_m\geq\frac{3}{2}m^2-\frac{7}{2}m+1$. In this case, the surface $X_m\subseteq\p^4$ is a counterexample of the Eisenbud-Goto conjecture.
	\end{thmalphabetintro}
	
	Here, $L(k)$ and $W(k)$ are the matrices whose rows correspond to some polynomials we consider. The invertibility of $L(k)$ implies the existence of a homogeneous polynomial containing the term $x_0^{\frac{3}{2}m^2-\frac{7}{2}m-1}x_1^2$ in $I_{X_m}$ (\Cref{lemma_x0tx1^2}) and the invertibility of $W(k)$ is equivalent to the nonexistence of homogeneous polynomials containing the term $x_0^{\frac{3}{2}m^2-\frac{7}{2}m-2}x_1^2$ in $I_{X_m}$ (\Cref{lemma_x0tx1^2_not}). We consider only when $m$ is divided by $6$ since the reason that $\maxdeg I_{X_m}$ is large can be explained uniformly under this condition as \Cref{mainthm_general_m_intro}.
	
	Geometrically, $X_m\subseteq\p^4$ is an outer projection of a toric variety of degree $m^2-m+3$ in $\p^5$ that is given by the image of the monomial rational map
	\[
	\varphi:\p^2=\proj\k[y_0,y_1,y_2] \dashrightarrow \p^5=\proj\k[x_0,\cdots,x_5]
	\] corresponding to $(\o_{\p^2}(m),y_1^{m-1}y_2,
	y_0^{m-1}y_1,
	y_1^m,y_1^{m-3}y_2^3,
	y_0^m,
	y_2^{m-1}y_0)$.
	
	Using \Cref{mainthm_general_m_intro}, we give a theoretical proof that $X_6$ is a counterexample to the regularity conjecture.
	
	\begin{thmalphabetintro}\label{mainthm_m=6_intro}
		The surface $X_6\subseteq\p^4$ is a counterexample to the Eisenbud-Goto conjecture as $\maxdeg X_6\geq 34$ and $\deg X_6=33$. It violates the normality conjecture as it is neither $31$-normal nor $32$-normal.
	\end{thmalphabetintro}
	
	It is hard to determine the regularity of $X_m$ for all $m\geq 6$. However we believe that
	\begin{equation*}
		\reg X_m= \begin{cases}
			\frac{3}{2}m^2-\frac{7}{2}m+1 &\text{if $3\mid m$}\\
			\frac{3}{2}m^2-\frac{7}{2}m+2 &\text{otherwise}
		\end{cases}
	\end{equation*}
	holds for all $m\geq 6$. We have checked this up to $m=25$ using Macaulay2 (\cite{M2}), and got the following result.
	{
		\footnotesize
		\[\begin{array}{|c|c|c|c|c|c|c|c|c|c|c|c|}
			\hline
			m        & 6  & 7  & 8  & 9  & 10  & \cdots & 21  & 22  & 23  & 24  & 25  \\ \hline
			\deg X_m & 33 & 45 & 59 & 75 & 93  & \cdots & 423 & 465 & 509 & 555 & 603 \\ \hline
			\reg X_m & 34 & 51 & 70 & 91 & 117 & \cdots & 589 & 651 & 715 & 781 & 852 \\ \hline
		\end{array}\]
	}
	
	\noindent
	Furthermore, we have checked the invertibility of $L(k)$ and $W(k)$ up to $k=27$ using Macaulay2. Hence we conclude the following.
	\vspace{-8pt}
	\begin{center}
		\resizebox{\linewidth}{!}{%
			$\begin{array}{|c|c|c|c|c|c|c|c|c|c|c|c|c|}
				\hline
				m        & 30        & 36        & 42        & 48        &  54    & \cdots & 138        & 144        & 150        & 156        & 162         \\ \hline
				\deg X_m & 873       & 1263      & 1725      & 2259      &  2865    & \cdots & 18909      & 20595      & 22353      & 24183      & 26085       \\ \hline
				\reg X_m & \geq 1246 & \geq 1819 & \geq 2500 & \geq 3289 & \geq 4186 & \cdots & \geq 28084 & \geq 30601 & \geq 33226 & \geq 35959 & \geq 38800 \\ \hline
			\end{array}$%
		}	
	\end{center}
	\vspace{3pt}
	
	Also, using Macaulay2, we get the following information of $X_6$:
	\begin{itemize}
		\item $d:=\deg X_6=33$, $e:=\codim X_6=2$
		\item $\maxdeg X_6=34>d$
		\item $\reg X_6=34>d-e+1$
		\item $\reg \o_X=18\leq d-e$
		\item $h^1(\mathcal{I}_{X_6}(31))=2$, $h^1(\mathcal{I}_{X_6}(32))=1$, and $h^1(\mathcal{I}_{X_6}(33))=0$
		\item The singular locus of $X_6$ is a union of three lines.
	\end{itemize}
	
	Note that one of the counterexamples suggested by J. McCullough and I. Peeva is a threefold of degree 31 and regularity 38 in $\p^5$, and one of the counterexamples given by J. Choe recently is a threefold of degree 18 and regularity 18 in $\p^5$.
	Using binomial rational maps between projective spaces, we give counterexamples of a surface and a threefold whose degrees and regularities are both 11, which are the counterexamples with the smallest regularity among those we are aware of (\Cref{example_reg11}, \Cref{example_reg11_threefold}).
	\\
	\\
	\indent
	In this paper, we work over an algebraically closed field $\k$ of characteristic zero unless stated otherwise and use the convention that varieties are integral. Likewise, curves, surfaces, and threefolds are assumed to be integral.
	For a projective variety $X\subseteq\p^r$, we denote its homogeneous coordinate ring as $S_X$ and its homogeneous ideal as $I_X$ unless stated otherwise.
	
	\section{Preliminaries}
	\noindent
	\textit{2.1.}
	\textbf{Known results and conjectures for toric varieties.}
	Despite of the counterexamples of the Eisenbud-Goto conjecture, it would be interesting to find the category of projective varieties in which the conjecture holds.
	In the next section, \Cref{mainthm_m=6} shows projective surfaces cannot be the candidate, hence the next varieties one may consider are projective toric varieties.
	
	Recall that the images of rational maps $\p^n\dashrightarrow\p^r$ given by linear systems $s_0,\cdots,s_r\in \H^0(\o_{\p^n}(d))$ where the $s_i$'s are all monomials are projective toric varieties.
	The family of surfaces we construct in the next section are the next simplest varieties to those in the following sense: those are images of rational maps between projective spaces given by linear systems whose sections are all monomials except one, which is a binomial.
	
	It has been widely believed that the Eisenbud-Goto conjecture is true for projective toric varieties although even the weak form of the conjecture is open currently. 
	
	\begin{conjecture}
		Let $X\subseteq\p^r$ be a projective toric variety. Then
		\begin{enumerate}
			\item (Weak inequality) $\maxdeg X\leq \deg X$;
			\item (Strong inequality) $\reg X\leq \deg X -\codim X +1$.
		\end{enumerate}
	\end{conjecture}
	
	In fact, for projective toric varieties of codimension two, the Eisenbud-Goto regularity conjecture holds as \Cref{thm_bound_lattice} shows.
	
	\begin{thm}[{\cite[Theorem 7.3]{PS}}]\label{thm_bound_lattice}
		Let $I$ be a homogeneous lattice ideal of a polynomial ring of codimension two that does not contain any linear forms. Then
		\[
		\reg I\leq \deg I,
		\]
		and the inequality holds strictly if $I$ is a toric ideal.
	\end{thm}
	~\\
	\noindent
	\textit{2.2.}
	\textbf{Partial elimination ideals.}
	In \cite{G}, M. Green introduced the concept of the partial elimination ideals which are deeply related to projections.
	As the counterexamples in the main theorems are images of projections of toric varieties from points, we use partial elimination ideals to observe the effect of the projection.
	For this purpose, we briefly review some properties of the partial elimination ideals.
	
	Let $R=\k[x_0,\cdots,x_r]$ and $S=\k[x_1,\cdots,x_r]$ be the polynomial rings.
	For a homogeneous ideal $I\subseteq R$, we define the set
	\[\footnotesize
	K_i(I):=\{f'\in S\mid \text{$\exists f\in I$ such that }f=f'x_0^i+\text{(lower $x_0$-degree terms)}\}\cup\{0\}.
	\]
	\noindent
	It is easy to see $K_i(I)$ is a homogeneous ideal of $S$, and it is called \textit{the $i$-th partial elimination ideal of $I$ with respect to $x_0$}.
	The following are some known properties of the partial elimination ideals.
	\begin{itemize}
		\item The filtration
		\[
		K_0(I)\subseteq K_1(I)\subseteq K_2(I)\subseteq\cdots\subseteq S
		\]
		stabilizes at some point so that $K_s(I)=K_{s+1}(I)=\cdots$. The minimum integer $s$ with $K_s(I)=K_{s+1}(I)=\cdots$ is called the stabilization number.
		\item Let $X\subseteq\p^r=\proj R$ be a projective variety and $q\in \p^r$ be a point. By a coordinate change, we may assume $q=[1,0,\cdots,0]$.
		Let
		\[\pi_q:X\dashrightarrow \p^{r-1}=\proj S\]
		be the projection of $X$ from $q$ and $X_q\subseteq\p^{r-1}$ be the image. Then we have $K_0(I_X)=I_{X_q}$, and $K_s(I_X)$ is an ideal that defines the projectivized tangent cone $\p TC_qX$ of $X$ at $q$ scheme-theoretically. Here, we use the convention $\p TC_qX$ is the empty set if $q\notin X$. Indeed, in that case, $K_s(I_X)=S$ and $s=\min\{\deg f\mid q\notin V(f),\ f\in I_X\}$.
		In any case, $K_i(I_X)$ defines \[\{p\in \p^{r-1}\mid \leng \pi_q^{-1}(p)>i\}\cup \p TC_qX\] set-theoretically.
	\end{itemize}
	
	Also, the following theorem shows the regularity of $X_q$ can be bounded using $\reg K_i(I_X)$.
	
	\begin{thm}[{\cite[Theorem 1.4]{KNV}}]\label{bound_reg_projection}
		Let $I\subsetneq R=\k[x_0,\cdots,x_r]$ be any nonzero homogeneous ideal.
		Then we have
		\[
		\reg_S K_0(I)\leq \max\limits_{1\leq i\leq s-1}(\reg_R I, \reg_S K_i(I)+i+1, \reg_S K_s(I)+s-1).
		\]
	\end{thm}
	
	To compute the partial elimination ideals in Macaulay2, one may use the following result.
	
	\begin{prop}\cite[Proposition 3.4]{CS}\label{PEI_by_GB}
		Let $G$ be a Gr\"{o}bner basis of a homogeneous ideal $I\subseteq R$ with respect to an elimination order that eliminates $x_0$.
		Then
		\[
		\{f'\in S\mid \textnormal{$\exists f\in G$ such that }f=f'x_0^j+\textnormal{(lower $x_0$-degree terms) where }j\leq i\}
		\]
		is a Gr\"{o}bner basis of $K_i(I)$.
	\end{prop}
	
	\section{Counterexamples to the Eisenbud-Goto conjecture}
	
	Let $X$ be a nondegenerate projective variety of degree $d$ and codimension $e$. The Eisenbud-Goto conjecture for $X$ is equivalent to $X$ satisfying both of the following condition:
	\begin{itemize}
		\item (normality conjecture) $X$ is $k$-normal for all $k\geq d-e$;
		\item (regularity conjecture for $\o_X$) $\reg\o_X\leq d-e$.
	\end{itemize}
	For smooth varieties, the second statement is proved by A. Noma, the second author, and J. Park (\cite{N2}, \cite{KP}).
	For surfaces, it is known that the second condition is always satisfied even if it is singular by the proposition below. Consequently, for surfaces, the Eisenbud-Goto regularity conjecture is equivalent to the normality conjecture.
	
	\begin{prop}[{\cite[Proposition 2.3]{NP1}}]\label{reg_of_surface}
		Let $X\subseteq\p^r$ be a nondegenerate surface of degree $d$ and codimension $e$ that is not necessarily smooth. Then $\reg \o_X\leq d-e$.
	\end{prop}
	\begin{proof}
		Let $C$ be the general hyperplane section of $X$. Then we have the exact sequence
		\[
		0\to \o_X(-1)\to\o_X\to\o_C\to 0.
		\]
		By tensoring with $\o_X(k+1)$ and taking the long exact sequence, we get
		{\footnotesize
			\[
			\cdots\to\H^{i-1}(C,\o_C(k+1))\to\H^i(X,\o_X(k))\to\H^i(X,\o_X(k+1))\to\H^i(C,\o_C(k+1))\to\cdots.
			\]
		}
		
		We claim $\H^i(X,\o_X(k))\hookrightarrow\H^i(X,\o_X(k+1))$ for all $i\geq 1$ and $k\geq d-e-i$.
		For the case $i\geq 2$, it follows easily as we have $\H^{i-1}(C,\o_C(k+1))=0$ for all $k\geq d-e-i$ since $\o_C$ is $(d-e)$-regular by \cite[Theorem 1.1]{GLP}.
		
		Now assume $i=1$. We have the following commutative diagram
		
		\[
		\begin{tikzcd}
			(S_X)_k\ar[r,"",hookrightarrow]\ar[d,"",twoheadrightarrow] & \H^0(X,\o_X(k))\ar[d,""] \\
			(S_X/L)_k\ar[r,""]        & \H^0(C,\o_C(k))
		\end{tikzcd}
		\]
		where $S_X$ denotes the homogeneous coordinate ring of $X\subseteq\p^r$ and $L\in S_X$ is the element of degree 1 corresponding to $\p^{r-1}$. For $k\geq d-e$, the map $(S_X/L)_k\to \H^0(C,\o_C(k))$ is surjective as $C$ is $k$-normal by \cite[Theorem 1.1]{GLP}. Then $\H^0(X,\o_X(k))\to\H^0(C,\o_C(k))$ is also surjective.
		Hence we get $\H^1(X,\o_X(k-1))\hookrightarrow\H^1(X,\o_X(k))$ for all $k\geq d-e$, and the claim holds also for $i=1$. Therefore the result follows by Serre's vanishing theorem.
	\end{proof}
	
	Now we prove the existence of a surface counterexample. For this purpose, we define some matrices:
	
	\begin{defn}\label{defn_N_L}
		Let $k\geq 1$ be an integer and
		\begin{align*}
			b_1(k) &= 6k^2-k-1\\
			b_2(k) &= 18k^2-9k\\
			b_3(k) &= 18k^2-10k+1\\
			b_4(k) &= 12k^2-7k.
		\end{align*}
		Then we define the following matrices:
		\begin{itemize}
			\item For $i=1,2,3,4$, define the $(b_i(k)+1)\times (54k^2-27k+3)$ matrix $N_i(k)$ as
			
			\NiceMatrixOptions{delimiters/max-width,code-for-first-row = \scriptscriptstyle}{
				\begin{center}
					\resizebox{30pt}{!}{%
						$N_1(k)=$%
					}\resizebox{\linewidth-40pt}{!}{%
						$\begin{bNiceMatrix}[first-row]
							1 						&2 						&\cdots					&
							&\cdots					&6k^2-k				&\cdots
							&						&						&
							&
							&\cdots					&18k^2-3k-2				&\cdots
							&54k^2-27k+3 															\\
							& 						&						&
							&						&1						&
							&						&						&
							&
							&						&						&
							&																		\\
							& 						&						&
							&1						&3						&3
							&1						&						&
							&
							&						&						&
							&																		\\
							& 						&						&1
							&6						&15						&20
							&15						&6						&1
							&
							&						&						&
							&																		\\
							&					 	&						&
							&						&\vdots					&
							&						&						&
							&
							&						&						&
							&																		\\
							1						&{3b_1(k)\choose 1}		&\cdots					&
							&						&						&
							&						&						&
							&\cdots
							&{3b_1(k)\choose 3b_1(k)-1}	&1						&						
							&						
						\end{bNiceMatrix}$%
					}
				\end{center}
				
				\begin{center}
					\resizebox{30pt}{!}{%
						$N_2(k)=$%
					}\resizebox{\linewidth-40pt}{!}{%
						$\begin{bNiceMatrix}[first-row]
							1 						&2 						&\cdots					&\cdots
							&18k^2-9k+1				&\cdots					&
							&						&						&
							&						&\cdots					&54k^2-27k+2
							&54k^2-27k+3															\\
							& 						&						&
							&1						&2						&1
							&						&						&
							&						&						&
							&																\\
							& 						&						&1
							&5						&10						&10
							&5						&1						&
							&						&						&
							&																\\
							&					 	&1						&8
							&28						&56						&70
							&56						&28						&8
							&1						&						&
							&																\\
							&					 	&						&
							&\vdots					&						&
							&						&						&
							&						&						&
							&																	\\
							1						&{3b_2(k)+2\choose 1}		&\cdots					&
							&						&						&
							&						&						&
							&						&\cdots					&{3b_2(k)+2\choose 3b_2(k)+1}
							&1						
						\end{bNiceMatrix}$%
					}
				\end{center}
				
				\begin{center}
					\resizebox{30pt}{!}{%
						$N_3(k)=$%
					}\resizebox{\linewidth-40pt}{!}{%
						$\begin{bNiceMatrix}[first-row]
							1 						&\cdots					&3k					&\cdots
							&						&\cdots					&18k^2-7k+1
							&\cdots					&						&
							&						&\cdots					&54k^2-27k+2
							&54k^2-27k+3															\\
							& 						&						&
							&						&						&1
							&						&						&
							&						&						&
							&																\\
							& 						&						&
							&						&1						&3
							&3						&1						&
							&						&						&
							&														\\
							&					 	&						&
							&1						&6						&15
							&20						&15						&6
							&1						&						&
							&														\\
							&					 	&						&
							&						&						&\vdots
							&						&						&
							&						&						&
							&																\\
							&						&1						&{3b_3(k)\choose 1}
							&\cdots					&						&
							&						&						&
							&						&\cdots					&{3b_3(k)\choose 3b_3(k)-1}
							&1						
						\end{bNiceMatrix}$%
					}
				\end{center}
				
				\begin{center}
					\resizebox{30pt}{!}{%
						$N_4(k)=$%
					}\resizebox{\linewidth-40pt}{!}{%
						$\begin{bNiceMatrix}[first-row]
							1 						&\cdots					&18k^2-6k+2				&\cdots
							&						&\cdots					&30k^2-13k+2
							&\cdots					&						&
							&						&						&\cdots
							&54k^2-27k+2			&54k^2-27k+3									\\
							& 						&						&
							&						&						&1
							&						&						&
							&						&						&
							&						&									\\
							& 						&						&
							&						&1						&3
							&3						&1						&
							&						&						&
							&						&									\\
							&					 	&						&
							&1						&6						&15
							&20						&15						&6
							&1						&						&
							&						&								\\
							&					 	&						&
							&						&						&\vdots
							&						&						&
							&						&						&
							&						&									\\
							&						&1						&{3b_4(k)\choose 1}
							&\cdots					&						&
							&						&						&
							&						&\cdots					&{3b_4(k)\choose 3b_4(k)-1}	&1						&						
						\end{bNiceMatrix}$%
					}
				\end{center}
			}
			\noindent
			where the numbers above the matrix indicate the indices of the columns.
			\item Define the $(54k^2-27k+4)\times (54k^2-27k+3)$ matrix \[
			N(k)=
			\begin{bmatrix}
				N_1(k)\\
				N_2(k)\\
				N_3(k)\\
				N_4(k)
			\end{bmatrix}_{\textstyle\raisebox{4pt}.}
			\]
			\item Define $L(k)$ to be the square matrix obtained by eliminating the first row of $N(k)$.
			\item For $1\leq i\leq 4$, define $W_i(k)$ to be the matrix obtained by eliminating the last row, the first column, and the last two columns from $N_i(k)$.
			\item Define the $(54k^2-27k)\times (54k^2-27k)$ matrix \[
			W(k)=
			\begin{bmatrix}
				W_1(k)\\
				W_2(k)\\
				W_3(k)\\
				W_4(k)
			\end{bmatrix}_{\textstyle\raisebox{4pt}.}
			\]
		\end{itemize}
	\end{defn}
	
	Then we give a family $\{X_m\}_{m\geq 6}$ of projective surfaces in $\p^4$ where we highly suspect $X_m$ is a counterexample of the regularity conjecture for all $m\geq 6$.
	
	\begin{thmalphabetmaintext}\label{mainthm_general_m}
		Let $m\geq 6$ be an integer and $X_m$ be the image of the binomial rational map \[\varphi_m:\p^2=\proj\k[y_0,y_1,y_2]\dashrightarrow\p^4=\proj\k[x_0,\cdots,x_4]\] given by the linear system \[(\o_{\p^2}(m),y_1^{m-1}y_2,
		y_0^{m-1}y_1,
		y_1^m+y_1^{m-3}y_2^3,
		y_0^m,
		y_2^{m-1}y_0).\]
		Then $X_m\subseteq\p^4$ is a surface of degree $m^2-m+3$.
		In addition, suppose $m=6k$ for some integer $k\geq 1$ and define the square matrices $L(k)$ and $W(k)$ as in \Cref{defn_N_L}. If $L(k)$ and $W(k)$ are both invertible, then $\maxdeg X_m\geq\frac{3}{2}m^2-\frac{7}{2}m+1$. In this case, the surface $X_m\subseteq\p^4$ is a counterexample of the Eisenbud-Goto conjecture.
	\end{thmalphabetmaintext}
	
	\begin{rem}\label{family_of_surfaces}
		Although there is a condition that requires the matrices $L(k)$ and $W(k)$ to be invertible, we believe it is satisfied for all $k\geq 1$.
		We have checked the invertibility of these matrices up to $k=27$ using Macaulay2.
	\end{rem}
	
	We found the family $\{X_m\}_{m\geq 6}$ by using Macaulay2 primarily.
	However, it is possible to prove theoretically that $X_6$ is a surface counterexample without depending on the integrity of computer algebra systems. Specifically, we prove the following using \Cref{mainthm_general_m}.
	
	\begin{thmalphabetmaintext}\label{mainthm_m=6}
		The surface $X_6\subseteq\p^4$ is a counterexample to the Eisenbud-Goto conjecture as $\maxdeg X_6\geq 34$ and $\deg X_6=33$. It violates the normality conjecture as it is neither $31$-normal nor $32$-normal.
	\end{thmalphabetmaintext}
	
	To show \Cref{mainthm_general_m} and \Cref{mainthm_m=6}, we give some lemmas and a definition.
	
	\begin{defn}
		Let $F=c_1M_1+\cdots+c_rM_r$ be a nonzero homogeneous polynomial in $\k[x_0,\cdots,x_n]$ where $M_i$'s are distinct monomials and $c_i$'s are nonzero scalars. Then we define \[\mon(F):=\{M_1,\cdots,M_r\}.\]
	\end{defn}
	
	\begin{lem}\label{lemma_x0t_x0tx1}
		Assume $m=6k$ and use the notations in \Cref{mainthm_general_m}. Then for any non-constant homogeneous polynomial $F\in I_{X_m}$, we have $x_0^t\notin \mon(F)$ and $x_0^tx_1\notin \mon(F)$ for all $t\geq 0$.
	\end{lem}
	\begin{proof}
		Let
		\begin{align*}
			s_0 &=y_1^{6k-1}y_2\\
			s_1 &=y_0^{6k-1}y_1\\
			s_2 &=y_1^{6k}+y_1^{6k-3}y_2^3=y_1^{6k-3}(y_1^3+y_2^3)\\
			s_3 &=y_0^{6k}\\
			s_4 &=y_2^{6k-1}y_0.
		\end{align*}
		As $I_{X_m}$ is the kernel of the map
		\begin{align*}
			\varphi^\#:\k[x_0,\cdots,x_4] &\to \k[y_0,y_1,y_2]\\
			x_i &\mapsto s_i,
		\end{align*}
		a homogeneous polynomial $F\in\k[x_0,\cdots,x_4]$ is contained in $I_{X_m}$ if and only if $F(s_0,\cdots,s_4)\in\k[y_0,y_1,y_2]$ is zero.
		
		Note that for a monomial $M\in\k[x_0,\cdots,x_4]$ of degree $t$, if $y_1^{a_1}y_2^{a_2}\in \mon(M(s_0,\cdots,s_4))$ for some integers $a_1,a_2$ such that $a_2\equiv t\pmod 3$, then $M=x_0^{t-3u}x_2^{3u}$ for some integer $0\leq u\leq \frac{t}{3}$. Furthermore, for such $M$, all elements of $\mon(M(s_0,\cdots,s_4))$ are of the form $y_1^{a_1}y_2^{a_2}$ where $a_2\equiv t\pmod 3$.
		
		Suppose there is a homogeneous polynomial $F\in I_{X_m}$ such that $x_0^t\in\mon(F)$. Since $F(s_0,\cdots,s_4)=0$, the term $s_0^t$ should be cancelled by other terms over $s_0,\cdots,s_4$. Also as $s_0^t=y_1^{(6k-1)t}y_2^t$, there should be a polynomial $G\in I_{X_m}$ such that $\mon(G)=\{x_0^t,M_1,\cdots,M_r\}$ where $M_i=x_0^{t-3u_i}x_2^{3u_i}$ for some integer $u_i\geq 1$ for all $i$. Indeed, such $G$ can be obtained by eliminating all monomials not of the form $x_0^{t-3u}x_2^{3u}$ from $F$.
		
		However, the elements of $U:=\{s_0^{t-3u}s_2^{3u}\mid 0\leq u\leq \left\lfloor\frac{t}{3}\right\rfloor=:d\}$ are $\k$-linearly independent as \[s_0^{t-3u}s_2^{3u}\vert_{y_1=1}=\sum_{j=0}^{3u}{3u\choose j}y_2^{t+6u-3j}\]
		so that $y_2$-degrees of the elements of $U$ are all different.
		Consequently, there is no homogeneous polynomial $F\in I_{X_m}$ such that $x_0^t\in\mon(F)$.
		
		Now we prove the second part using the similar argument. 
		Note that for a monomial $M\in\k[x_0,\cdots,x_4]$ of degree $t+1$, if $y_0^{6k-1}y_1^{a_1}y_2^{a_2}\in \mon(M(s_0,\cdots,s_4))$ for some integers $a_1,a_2$ such that $a_2\equiv t\pmod 3$, then $M=x_0^{t-3u}x_1x_2^{3u}$ or $M=x_0^{t-6k+2-3u}x_2^{3u}x_4^{6k-1}$ for some integer $u\geq 0$.
		Furthermore, for such $M$, all elements of $\mon(M(s_0,\cdots,s_4))$ are of the form $y_0^{6k-1}y_1^{a_1}y_2^{a_2}$ where $a_2\equiv t\pmod 3$.
		
		Suppose there is a homogeneous polynomial $F\in I_{X_m}$ such that $x_0^tx_1\in\mon(F)$. Since $F(s_0,\cdots,s_4)=0$, the term $s_0^ts_1$ should be cancelled by other terms over $s_0,\cdots,s_4$. As $s_0^ts_1=y_0^{6k-1}y_1^{(6k-1)t+1}y_2^t$, there should be a polynomial $G\in I_{X_m}$ such that $\mon(G)=\{x_0^tx_1,M_1,\cdots,M_r\}$ where $M_i=x_0^{t-3u_i}x_1x_2^{3u_i}$ for some $u_i\geq 1$ or $M_i=x_0^{t-6k+2-3u_i}x_2^{3u_i}x_4^{6k-1}$ for some $u_i\geq 0$.
		
		We claim $s_0^ts_1,M_1(s_0,\cdots,s_4),\cdots,M_r(s_0,\cdots,s_4)\in\k[y_0,y_1,y_2]$ are $\k$-linearly independent to induce a contradiction.
		Define $d_1:=\left\lfloor\frac{t}{3}\right\rfloor$, $d_2:=\left\lfloor\frac{t-6k+2}{3}\right\rfloor$, and \[T:=\{s_0^{t-3u}s_1s_2^{3u}\mid 0\leq u\leq d_1\}\cup\{s_0^{t-6k+2-3u}s_2^{3u}s_4^{6k-1}\mid 0\leq u\leq d_2\}.\]
		Then it is enough to show the elements of $T$ are $\k$-linearly independent.
		For $0\leq u\leq d_1$, as
		\[
			s_0^{t-3u}s_1s_2^{3u}\vert_{y_0=1,y_1=1} =\sum_{j=0}^{3u}{3u\choose j}y_2^{t-3u+3j},
		\]
		the $y_2$-degree of $s_0^{t-3u}s_1s_2^{3u}$ is $t+6u$. Also for $0\leq u\leq d_2$, as
		\[
			s_0^{t-6k+2-3u}s_2^{3u}s_4^{6k-1}\vert_{y_0=1,y_1=1} =\sum_{j=0}^{3u}{3u\choose j}y_2^{3c+t-3u+3j},
		\]
		the $y_2$-degree of $s_0^{t-6k+2-3u}s_2^{3u}s_4^{6k-1}$ is $t+6u+3c$ where $c:=12k^2-6k+1$ is an odd number.
		Therefore, $y_2$-degrees of the elements of $T$ are all different so that those are $\k$-linearly independent.
		Consequently, there is no homogeneous polynomial $F\in I_{X_m}$ such that $x_0^tx_1\in\mon(F)$.
	\end{proof}
	
	\begin{lem}\label{lemma_x0tx1^2}
		Assume $m=6k$ and use the notations in \Cref{mainthm_general_m}. If $L(k)$ is invertible, then there exists a homogeneous polynomial $F\in I_{X_m}$ such that $x_0^{\frac{3}{2}m^2-\frac{7}{2}m-1}x_1^2\in \mon(F)$.
	\end{lem}
	\begin{proof}
		Suppose $L(k)$ is invertible, and note that \[s_0^{\frac{3}{2}m^2-\frac{7}{2}m-1}s_1^2=y_0^{12k-2}y_1^{(6k-1)(54k^2-21k-1)+2}y_2^{54k^2-21k-1}.\]
		For a monomial $M\in\k[x_0,\cdots,x_4]$ of degree $\frac{3}{2}m^2-\frac{7}{2}m+1=54k^2-21k+1$, if $y_0^{12k-2}y_1^{a_1}y_2^{a_2}\in \mon(M(s_0,\cdots,s_4))$ for some integers $a_1,a_2$ such that $a_2\equiv 2\pmod 3$, then $M$ must be one of the following types:
		\begin{enumerate}
			\item $M=x_0^{54k^2-21k-3u-1}x_1^2x_2^{3u}$ for some $0\leq u\leq 18k^2-7k-1$;
			\item $M=x_0^{54k^2-27k-3u}x_2^{3u+2}x_3x_4^{6k-2}$ for some $0\leq u\leq 18k^2-9k$;
			\item $M=x_0^{54k^2-27k-3u+1}x_1x_2^{3u}x_4^{6k-1}$ for some $0\leq u\leq 18k^2-9k$; or
			\item $M=x_0^{54k^2-33k-3u+3}x_2^{3u}x_4^{12k-2}$ for some $0\leq u\leq 18k^2-11k+1$.
		\end{enumerate}
		Furthermore, if $M$ is one of the types above, then the elements of $\mon(M(s_0,\cdots,s_4))$ are of the form $y_0^{12k-2}y_1^{a_1}y_2^{a_2}$ where $a_2\equiv 2\pmod 3$.
		Define
		\begin{align*}
			T_1 &:= \{s_0^{54k^2-21k-3u-1}s_1^2s_2^{3u}\mid 0\leq u\leq b_1(k)=6k^2-k-1\} \\
			T_2 &:= \{s_0^{54k^2-27k-3u}s_2^{3u+2}s_3s_4^{6k-2}\mid 0\leq u\leq b_2(k)=18k^2-9k\} \\
			T_3 &:= \{s_0^{54k^2-27k-3u+1}s_1s_2^{3u}s_4^{6k-1}\mid 0\leq u\leq b_3(k)=18k^2-10k+1\} \\
			T_4 &:= \{s_0^{54k^2-33k-3u+3}s_2^{3u}s_4^{12k-2}\mid 0\leq u\leq b_4(k)=12k^2-7k\}.
		\end{align*}
		
		To show there is a homogeneous polynomial $F\in I_{X_m}$ such that $x_0^{54k^2-21k-1}x_1^2\in\mon(F)$, it is enough to show the elements of $T:=T_1\cup T_2\cup T_3\cup T_4$ are $\k$-linearly dependent while the elements of $T\backslash\{s_0^{54k^2-21k-1}s_1^2\}$ are $\k$-linearly independent.
		As every element in $T$ has the same degree and also the same $y_0$-degree, it is enough to prove it after taking the substitution $y_0=y_1=1$ to every elements of $T$.
		Define
		\begin{align*}
			f_1(u) :=& s_0^{54k^2-21k-3u-1}s_1^2s_2^{3u}\vert_{y_0=1,y_1=1}\\
			=& \sum_{j=0}^{3u}{3u\choose j}y_2^{54k^2-21k-1-3u+3j}\text{ (for $0\leq u\leq b_1(k)$)}\\[6pt]
			f_2(u) :=& s_0^{54k^2-27k-3u}s_2^{3u+2}s_3s_4^{6k-2}\vert_{y_0=1,y_1=1}\\
			=& \sum_{j=0}^{3u+2}{3u+2\choose j}y_2^{90k^2-45k+2-3u+3j}\text{ (for $0\leq u\leq b_2(k)$)}\\[6pt]
			f_3(u) :=& s_0^{54k^2-27k-3u+1}s_1s_2^{3u}s_4^{6k-1}\vert_{y_0=1,y_1=1}\\
			=& \sum_{j=0}^{3u}{3u\choose j}y_2^{90k^2-39k+2-3u+3j}\text{ (for $0\leq u\leq b_3(k)$)}\\[6pt]
			f_4(u) :=& s_0^{54k^2-33k-3u+3}s_2^{3u}s_4^{12k-2}\vert_{y_0=1,y_1=1}\\
			=& \sum_{j=0}^{3u}{3u\choose j}y_2^{126k^2-57k+5-3u+3j}\text{ (for $0\leq u\leq b_4(k)$)}.
		\end{align*}
		
		Note that the lowest $y_2$-degree of elements in $T$ is $36k^2-18k+2$ and the highest $y_2$-degree of elements in $T$ is $198k^2-99k+8$.
		Let $\widetilde{N}_i(k)$ be a $(b_i(k)+1)\times (198k^2-99k+9)$ matrix whose $(u,v)$-th entry is the coefficient of $y_2^{v-1}$ in $f_i(u-1)$ for $i=1,2,3,4$.
		As the $j$-th column of $\widetilde{N}_i(k)$ is zero whenever $j\not\equiv 0 \pmod 3$ or $1\leq j\leq 36k^2-18k+2$, we remove all such columns from $\widetilde{N}_i(k)$. By this elimination, we get the matrices $N_1(k),N_2(k),N_3(k),N_4(k)$. As the rows of $N(k)$ represent the coefficients of the elements of $T$ and the rows of $L(k)$ represent the coefficients of the elements of $T\backslash\{s_0^{54k^2-21k-1}s_1^2\}$, we get the result as $L(k)$ is invertible and the number of rows of $N(k)$ is larger than the number of columns of $N(k)$.
	\end{proof}
	
	\begin{lem}\label{lemma_x0tx1^2_not}
		Assume $m=6k$ and use the notations in \Cref{mainthm_general_m}. Then $W(k)$ is invertible if and only if $x_0^{\frac{3}{2}m^2-\frac{7}{2}m-2}x_1^2\notin \mon(F)$ for any homogeneous polynomial $F\in I_{X_m}$.
	\end{lem}
	\begin{proof}
		We argue as before. Note that \[s_0^{\frac{3}{2}m^2-\frac{7}{2}m-2}s_1^2=y_0^{12k-2}y_1^{(6k-1)(54k^2-21k-2)+2}y_2^{54k^2-21k-2}.\]
		For a monomial $M\in\k[x_0,\cdots,x_4]$ of degree $\frac{3}{2}m^2-\frac{7}{2}m=54k^2-21k$, if $y_0^{12k-2}y_1^{a_1}y_2^{a_2}\in \mon(M(s_0,\cdots,s_4))$ for some integers $a_1$ and $a_2$ with $a_2\equiv 1\pmod 3$, then $M$ must be one of the following types:
		\begin{enumerate}[leftmargin=1cm]
			\item $M=x_0^{54k^2-21k-3u-2}x_1^2x_2^{3u}$ for some $0\leq u\leq 18k^2-7k-1$;
			\item $M=x_0^{54k^2-27k-3u-1}x_2^{3u+2}x_3x_4^{6k-2}$ for some $0\leq u\leq 18k^2-9k-1$;
			\item $M=x_0^{54k^2-27k-3u}x_1x_2^{3u}x_4^{6k-1}$ for some $0\leq u\leq 18k^2-9k$; or
			\item $M=x_0^{54k^2-33k-3u+2}x_2^{3u}x_4^{12k-2}$ for some $0\leq u\leq 18k^2-11k$.
		\end{enumerate}
		Furthermore, if $M$ is one of the types above, then the elements of $\mon(M(s_0,\cdots,s_4))$ are of the form $y_0^{12k-2}y_1^{a_1}y_2^{a_2}$ where $a_2\equiv 1\pmod 3$.
		Define the sets
		\begin{align*}
			T_1 &:= \{s_0^{54k^2-21k-3u-2}s_1^2s_2^{3u}\mid 0\leq u\leq 18k^2-7k-1\} \\
			T_2 &:= \{s_0^{54k^2-27k-3u-1}s_2^{3u+2}s_3s_4^{6k-2}\mid 0\leq u\leq 18k^2-9k-1\} \\
			T_3 &:= \{s_0^{54k^2-27k-3u}s_1s_2^{3u}s_4^{6k-1}\mid 0\leq u\leq 18k^2-9k\} \\
			T_4 &:= \{s_0^{54k^2-33k-3u+2}s_2^{3u}s_4^{12k-2}\mid 0\leq u\leq 18k^2-11k\}.
		\end{align*}
		
		To see whether there is no homogeneous polynomial $F\in I_{X_m}$ such that $x_0^{54k^2-21k-1}x_1^2\in\mon(F)$, it is enough to see whether the elements of $T:=T_1\cup T_2\cup T_3\cup T_4$ are $\k$-linearly independent.
		As every element in $T$ has the same degree and also the same $y_0$-degree, it is enough to check the independence after taking the substitution $y_0=y_1=1$ to every elements of $T$.
		Define
		\begin{align*}
			f_1(u) :=& s_0^{54k^2-21k-3u-2}s_1^2s_2^{3u}\vert_{y_0=1,y_1=1}\\
			=& \sum_{j=0}^{3u}{3u\choose j}y_2^{54k^2-21k-2-3u+3j}\text{ (for $0\leq u\leq 18k^2-7k-1$)}\\[6pt]
			f_2(u) :=& s_0^{54k^2-27k-3u-1}s_2^{3u+2}s_3s_4^{6k-2}\vert_{y_0=1,y_1=1}\\
			=& \sum_{j=0}^{3u+2}{3u+2\choose j}y_2^{90k^2-45k+1-3u+3j}\text{ (for $0\leq u\leq 18k^2-9k-1$)}\\[6pt]
			f_3(u) :=& s_0^{54k^2-27k-3u}s_1s_2^{3u}s_4^{6k-1}\vert_{y_0=1,y_1=1}\\
			=& \sum_{j=0}^{3u}{3u\choose j}y_2^{90k^2-39k+1-3u+3j}\text{ (for $0\leq u\leq 18k^2-9k$)}\\[6pt]
			f_4(u) :=& s_0^{54k^2-33k-3u+2}s_2^{3u}s_4^{12k-2}
			\vert_{y_0=1,y_1=1}\\
			=& \sum_{j=0}^{3u}{3u\choose j}y_2^{126k^2-57k+4-3u+3j}\text{ (for $0\leq u\leq 18k^2-11k$)}.
		\end{align*}
		
		Let $T'=\{s\vert_{y_0=1,y_1=1}\mid s\in T\}$.
		To see whether the elements of $T'$ are linearly independent, we do not need to consider all elements of $T'$. For instance, $f_1(18k^2-7k-1)$ is the only element in $T'$ that contains the term ${y_2}^1$ so that any $\k$-linear combination of elements of $T'$ containing $f_1(18k^2-7k-1)$ cannot be zero.
		After repeating the elimination of such elements, it is enough to see whether the elements of the following set are linearly independent:
		\begin{alignat*}{3}
			&\{f_1(u)\mid 0\leq u\leq 6k^2-k-2\}\ &&\cup \ &&\{f_2(u)\mid 0\leq u\leq 18k^2-9k-1\} \\
			\cup\ &\{f_3(u)\mid 0\leq u\leq 18k^2-10k\}\ &&\cup \ &&\{f_4(u)\mid 0\leq u\leq 12k^2-7k-1\}.
		\end{alignat*}
		We get the matrix $W(k)$ by constructing the coefficient matrix as before. Hence the result follows.
	\end{proof}
	
	Now we prove the main theorems.
	
	\begin{proof}[Proof of \Cref{mainthm_general_m}]
		Let $Y_m$ be the image of the rational map
		\[
		\p^2=\proj\k[y_0,y_1,y_2]\dashrightarrow\p^5=\proj\k[x_0,\cdots,x_5]
		\]
		given by $[y_1^{m-1}y_2,
		y_0^{m-1}y_1,
		y_1^m,y_1^{m-3}y_2^3,
		y_0^m,
		y_2^{m-1}y_0]$.
		Then $Y_m$ is a projective toric variety corresponding to \[\mathcal{A}:=\{(0,1),(m-1,0),(0,0),(0,3),(m,0),(1,m-1)\}\subseteq\z^2.\]
		Note that the lattice affinely generated by $\mathcal{A}$ is $\z^2$ and the normalized volume of the polytope $\conv(\mathcal{A})$ is $2!\frac{m^2-m+3}{2}=m^2-m+3$. Thus $\deg Y_m=m^2-m+3$ by Kushnirenko's theorem (\cite{K1}, also \textit{c.f.} \cite[Theorem 3.16]{T}).
		As $X_m$ is a birational outer projection of $Y_m$, we get $\deg X_m=m^2-m+3$.
		
		Now assume $m=6k$ for some integer $k\geq 1$.
		By \Cref{lemma_x0tx1^2}, there is a homogeneous polynomial $F\in I_{X_m}$ such that $x_0^{\frac{3}{2}m^2-\frac{7}{2}m-1}x_1^2\in \mon(F)$. On the other hand, for any non-constant polynomial $G\in I_{X_m}$ and for any $t\geq 0$, we have $x_0^t\notin \mon(G)$, $x_0^tx_1\notin \mon(G)$ by \Cref{lemma_x0t_x0tx1}, and $x_0^{\frac{3}{2}m^2-\frac{7}{2}m-2}x_1^2\notin \mon(G)$ by \Cref{lemma_x0tx1^2_not}.
		This shows that in any minimal generating set of $I_{X_m}$, there exists a homogeneous polynomial $F$ such that $x_0^{\frac{3}{2}m^2-\frac{7}{2}m-1}x_1^2\in \mon(F)$.
		Hence $\maxdeg I_{X_m}\geq\frac{3}{2}m^2-\frac{7}{2}m+1$.
	\end{proof}
	
	\begin{proof}[Proof of \Cref{mainthm_m=6}]
		Note that the entries of $L(1)$ and $W(1)$ are contained in $\z$. Let $L(1)_{\z/3\z}$ and $W(1)_{\z/3\z}$ be the matrices over $\z/3\z$ obtained by taking the quotient to every entry of $L(1)$ and $W(1)$.
		Then it is enough to show $L(1)_{\z/3\z}$ and $W(1)_{\z/3\z}$ are invertible over $\z/3\z$.
		
		We define $N_i(1)_{\z/3\z}$ in the same way above for $1\leq i\leq 4$.
		To calculate $L(1)_{\z/3\z}$ and $W(1)_{\z/3\z}$, note that it is enough to calculate the entries of $N_2(1)_{\z/3\z}$ and $N_3(1)_{\z/3\z}$ since the entries of all other matrices are induced from those.
		For $N_3(1)_{\z/3\z}$, note that $\binom{3a}{b}\equiv 0\pmod 3$ if and only if
		{\small
		\[
		\left\lfloor\frac{3a}{3}\right\rfloor+\left\lfloor\frac{3a}{3^2}\right\rfloor+\cdots
		>
		\left\lfloor\frac{b}{3}\right\rfloor+\left\lfloor\frac{b}{3^2}\right\rfloor+\cdots+
		\left\lfloor\frac{3a-b}{3}\right\rfloor+\left\lfloor\frac{3a-b}{3^2}\right\rfloor+\cdots
		\]
		}
		which is always the case whenever
		\begin{itemize}
			\item $3\nmid b$, or
			\item $3\mid a$ and $9\nmid b$, or
			\item $9\mid a$ and $27\nmid b$.
		\end{itemize}
		Also for $N_2(1)_{\z/3\z}$, note that $\binom{3a+2}{b}\equiv 0\pmod 3$ if and only if
		{\Small
		\[
		\left\lfloor\frac{3a+2}{3}\right\rfloor+\left\lfloor\frac{3a+2}{3^2}\right\rfloor+\cdots
		>
		\left\lfloor\frac{b}{3}\right\rfloor+\left\lfloor\frac{b}{3^2}\right\rfloor+\cdots+
		\left\lfloor\frac{3a+2-b}{3}\right\rfloor+\left\lfloor\frac{3a+2-b}{3^2}\right\rfloor+\cdots		
		\]
		}
		which is always the case whenever
		\begin{itemize}
			\item $3a+2\geq 9$ and $b<9$ and $3a+2-b<9$, or
			\item $3a+2\geq 18$ and $b<9$ and $3a+2-b<18$, or
			\item $3a+2\geq 27$ and $b<27$ and $3a+2-b<27$.
		\end{itemize}
		
		Using these and symmetry reduces the calculation a lot. After the calculations, we get the matrix $L(1)_{\z/3\z}$
		\begin{center}
			\resizebox{\linewidth}{!}{%
				\NiceMatrixOptions{rules/color=[gray]{0.75}, hvlines}
				$\begin{bNiceMatrix}[first-row,first-col]
					\CodeBefore
					\rectangleEntryColor[1-5][1-7]
					\rectangleEntryColor[2-4][2-9]
					\rectangleEntryColor[3-3][3-11]
					\rectangleEntryColor[4-2][4-13]
					\rectangleEntryColor[5-10][5-11]
					\rectangleEntryColor[6-9][6-13]
					\rectangleEntryColor[7-8][7-15]
					\rectangleEntryColor[8-7][8-17]
					\rectangleEntryColor[9-6][9-19]
					\rectangleEntryColor[10-5][10-21]
					\rectangleEntryColor[11-4][11-23]
					\rectangleEntryColor[12-3][12-25]
					\rectangleEntryColor[13-2][13-27]
					\rectangleEntryColor[14-1][14-29]
					\rectangleEntryColor[16-11][16-13]
					\rectangleEntryColor[17-10][17-15]
					\rectangleEntryColor[18-9][18-17]
					\rectangleEntryColor[19-8][19-19]
					\rectangleEntryColor[20-7][20-21]
					\rectangleEntryColor[21-6][21-23]
					\rectangleEntryColor[22-5][22-25]
					\rectangleEntryColor[23-4][23-27]
					\rectangleEntryColor[24-3][24-29]
					\rectangleEntryColor[25-19][25-19]
					\rectangleEntryColor[26-18][26-21]
					\rectangleEntryColor[27-17][27-23]
					\rectangleEntryColor[28-16][28-25]
					\rectangleEntryColor[29-15][29-27]
					\rectangleEntryColor[30-14][30-29]
					\cellcolor[HTML]{CCCCCC}{1-4,2-3,3-2,4-1,5-12,6-14,7-16,8-18,9-20,10-22,11-24,12-26,13-28,14-30,15-12,16-14,17-16,18-18,19-20,20-22,21-24,22-26,23-28,24-30}
					\Body
					& \matinds[01]          & \matinds[02]          & \matinds[03]          & \matinds[04]  & \matinds[05]  & \matinds[06]  & \matinds[07]  & \matinds[08]  & \matinds[09]  & \matinds[10] & \matinds[11]         & \matinds[12] & \matinds[13]         & \matinds[14] & \matinds[15]         & \matinds[16] & \matinds[17]         & \matinds[18] & \matinds[19]         & \matinds[20] & \matinds[21]         & \matinds[22] & \matinds[23]         & \matinds[24] & \matinds[25]         & \matinds[26] & \matinds[27]         & \matinds[28] & \matinds[29]         & \matinds[30]          \\
					\matinds[01]  & \entryzero          & \entryzero          & \entryzero          & \boldone & \entryzero  & \entryzero  & \coloredone  & \entryzero  & \entryzero  & \entryzero  & \entryzero  & \entryzero          & \entryzero  & \entryzero          & \entryzero  & \entryzero          & \entryzero  & \entryzero          & \entryzero  & \entryzero          & \entryzero  & \entryzero          & \entryzero  & \entryzero          & \entryzero  & \entryzero          & \entryzero  & \entryzero          & \entryzero  & \entryzero          \\
					\matinds[02]  & \entryzero          & \entryzero          & \boldone & \entryzero          & \entryzero  & \coloredtwo & \entryzero  & \entryzero  & \coloredone  & \entryzero  & \entryzero  & \entryzero          & \entryzero  & \entryzero          & \entryzero  & \entryzero          & \entryzero  & \entryzero          & \entryzero  & \entryzero          & \entryzero  & \entryzero          & \entryzero  & \entryzero          & \entryzero  & \entryzero          & \entryzero  & \entryzero          & \entryzero  & \entryzero          \\
					\matinds[03]  & \entryzero          & \boldone & \entryzero          & \entryzero          & \entryzero  & \entryzero  & \entryzero  & \entryzero  & \entryzero  & \entryzero  & \coloredone  & \entryzero          & \entryzero  & \entryzero          & \entryzero  & \entryzero          & \entryzero  & \entryzero          & \entryzero  & \entryzero          & \entryzero  & \entryzero          & \entryzero  & \entryzero          & \entryzero  & \entryzero          & \entryzero  & \entryzero          & \entryzero  & \entryzero          \\
					\matinds[04]  & \boldone & \entryzero          & \entryzero          & \coloredone          & \entryzero  & \entryzero  & \entryzero  & \entryzero  & \entryzero  & \coloredone  & \entryzero  & \entryzero          & \coloredone  & \entryzero          & \entryzero  & \entryzero          & \entryzero  & \entryzero          & \entryzero  & \entryzero          & \entryzero  & \entryzero          & \entryzero  & \entryzero          & \entryzero  & \entryzero          & \entryzero  & \entryzero          & \entryzero  & \entryzero          \\
					\matinds[05]  & \entryzero          & \entryzero          & \entryzero          & \entryzero          & \entryzero  & \entryzero  & \entryzero  & \entryzero  & \entryzero  & \coloredone  & \coloredtwo & \boldone & \entryzero  & \entryzero          & \entryzero  & \entryzero          & \entryzero  & \entryzero          & \entryzero  & \entryzero          & \entryzero  & \entryzero          & \entryzero  & \entryzero          & \entryzero  & \entryzero          & \entryzero  & \entryzero          & \entryzero  & \entryzero          \\
					\matinds[06]  & \entryzero          & \entryzero          & \entryzero          & \entryzero          & \entryzero  & \entryzero  & \entryzero  & \entryzero  & \coloredone  & \coloredtwo & \coloredone  & \coloredone          & \coloredtwo & \boldone & \entryzero  & \entryzero          & \entryzero  & \entryzero          & \entryzero  & \entryzero          & \entryzero  & \entryzero          & \entryzero  & \entryzero          & \entryzero  & \entryzero          & \entryzero  & \entryzero          & \entryzero  & \entryzero          \\
					\matinds[07]  & \entryzero          & \entryzero          & \entryzero          & \entryzero          & \entryzero  & \entryzero  & \entryzero  & \coloredone  & \coloredtwo & \coloredone  & \coloredtwo & \coloredone          & \coloredtwo & \coloredone          & \coloredtwo & \boldone & \entryzero  & \entryzero          & \entryzero  & \entryzero          & \entryzero  & \entryzero          & \entryzero  & \entryzero          & \entryzero  & \entryzero          & \entryzero  & \entryzero          & \entryzero  & \entryzero          \\
					\matinds[08]  & \entryzero          & \entryzero          & \entryzero          & \entryzero          & \entryzero  & \entryzero  & \coloredone  & \coloredtwo & \coloredone  & \entryzero  & \entryzero  & \entryzero          & \entryzero  & \entryzero          & \entryzero  & \coloredone          & \coloredtwo & \boldone & \entryzero  & \entryzero          & \entryzero  & \entryzero          & \entryzero  & \entryzero          & \entryzero  & \entryzero          & \entryzero  & \entryzero          & \entryzero  & \entryzero          \\
					\matinds[09]  & \entryzero          & \entryzero          & \entryzero          & \entryzero          & \entryzero  & \coloredone  & \coloredtwo & \coloredone  & \coloredone  & \coloredtwo & \coloredone  & \entryzero          & \entryzero  & \entryzero          & \coloredone  & \coloredtwo         & \coloredone  & \coloredone          & \coloredtwo & \boldone & \entryzero  & \entryzero          & \entryzero  & \entryzero          & \entryzero  & \entryzero          & \entryzero  & \entryzero          & \entryzero  & \entryzero          \\
					\matinds[10]  & \entryzero          & \entryzero          & \entryzero          & \entryzero          & \coloredone  & \coloredtwo & \coloredone  & \coloredtwo & \coloredone  & \coloredtwo & \coloredone  & \coloredtwo         & \coloredone  & \coloredone          & \coloredtwo & \coloredone          & \coloredtwo & \coloredone          & \coloredtwo & \coloredone          & \coloredtwo & \boldone & \entryzero  & \entryzero          & \entryzero  & \entryzero          & \entryzero  & \entryzero          & \entryzero  & \entryzero          \\
					\matinds[11] & \entryzero          & \entryzero          & \entryzero          & \coloredone          & \coloredtwo & \coloredone  & \entryzero  & \entryzero  & \entryzero  & \entryzero  & \entryzero  & \entryzero          & \coloredtwo & \coloredone          & \coloredtwo & \entryzero          & \entryzero  & \entryzero          & \entryzero  & \entryzero          & \entryzero  & \coloredone          & \coloredtwo & \boldone & \entryzero  & \entryzero          & \entryzero  & \entryzero          & \entryzero  & \entryzero          \\
					\matinds[12] & \entryzero          & \entryzero          & \coloredone          & \coloredtwo         & \coloredone  & \coloredone  & \coloredtwo & \coloredone  & \entryzero  & \entryzero  & \entryzero  & \coloredtwo         & \coloredone  & \coloredtwo         & \coloredtwo & \coloredone          & \coloredtwo & \entryzero          & \entryzero  & \entryzero          & \coloredone  & \coloredtwo         & \coloredone  & \coloredone          & \coloredtwo & \boldone & \entryzero  & \entryzero          & \entryzero  & \entryzero          \\
					\matinds[13] & \entryzero          & \coloredone          & \coloredtwo         & \coloredone          & \coloredtwo & \coloredone  & \coloredtwo & \coloredone  & \coloredtwo & \coloredone  & \coloredtwo & \coloredone          & \coloredtwo & \coloredone          & \coloredtwo & \coloredone          & \coloredtwo & \coloredone          & \coloredtwo & \coloredone          & \coloredtwo & \coloredone          & \coloredtwo & \coloredone          & \coloredtwo & \coloredone          & \coloredtwo & \boldone & \entryzero  & \entryzero          \\
					\matinds[14] & \coloredone          & \coloredtwo         & \coloredone          & \entryzero          & \entryzero  & \entryzero  & \entryzero  & \entryzero  & \entryzero  & \entryzero  & \entryzero  & \entryzero          & \entryzero  & \entryzero          & \entryzero  & \entryzero          & \entryzero  & \entryzero          & \entryzero  & \entryzero          & \entryzero  & \entryzero          & \entryzero  & \entryzero          & \entryzero  & \entryzero          & \entryzero  & \coloredone          & \coloredtwo & \boldone \\
					\matinds[15] & \entryzero          & \entryzero          & \entryzero          & \entryzero          & \entryzero  & \entryzero  & \entryzero  & \entryzero  & \entryzero  & \entryzero  & \entryzero  & \boldone & \entryzero  & \entryzero          & \entryzero  & \entryzero          & \entryzero  & \entryzero          & \entryzero  & \entryzero          & \entryzero  & \entryzero          & \entryzero  & \entryzero          & \entryzero  & \entryzero          & \entryzero  & \entryzero          & \entryzero  & \entryzero          \\
					\matinds[16] & \entryzero          & \entryzero          & \entryzero          & \entryzero          & \entryzero  & \entryzero  & \entryzero  & \entryzero  & \entryzero  & \entryzero  & \coloredone  & \entryzero          & \entryzero  & \boldone & \entryzero  & \entryzero          & \entryzero  & \entryzero          & \entryzero  & \entryzero          & \entryzero  & \entryzero          & \entryzero  & \entryzero          & \entryzero  & \entryzero          & \entryzero  & \entryzero          & \entryzero  & \entryzero          \\
					\matinds[17] & \entryzero          & \entryzero          & \entryzero          & \entryzero          & \entryzero  & \entryzero  & \entryzero  & \entryzero  & \entryzero  & \coloredone  & \entryzero  & \entryzero          & \coloredtwo & \entryzero          & \entryzero  & \boldone & \entryzero  & \entryzero          & \entryzero  & \entryzero          & \entryzero  & \entryzero          & \entryzero  & \entryzero          & \entryzero  & \entryzero          & \entryzero  & \entryzero          & \entryzero  & \entryzero          \\
					\matinds[18] & \entryzero          & \entryzero          & \entryzero          & \entryzero          & \entryzero  & \entryzero  & \entryzero  & \entryzero  & \coloredone  & \entryzero  & \entryzero  & \entryzero          & \entryzero  & \entryzero          & \entryzero  & \entryzero          & \entryzero  & \boldone & \entryzero  & \entryzero          & \entryzero  & \entryzero          & \entryzero  & \entryzero          & \entryzero  & \entryzero          & \entryzero  & \entryzero          & \entryzero  & \entryzero          \\
					\matinds[19] & \entryzero          & \entryzero          & \entryzero          & \entryzero          & \entryzero  & \entryzero  & \entryzero  & \coloredone  & \entryzero  & \entryzero  & \coloredone  & \entryzero          & \entryzero  & \entryzero          & \entryzero  & \entryzero          & \coloredone  & \entryzero          & \entryzero  & \boldone & \entryzero  & \entryzero          & \entryzero  & \entryzero          & \entryzero  & \entryzero          & \entryzero  & \entryzero          & \entryzero  & \entryzero          \\
					\matinds[20] & \entryzero          & \entryzero          & \entryzero          & \entryzero          & \entryzero  & \entryzero  & \coloredone  & \entryzero  & \entryzero  & \coloredtwo & \entryzero  & \entryzero          & \coloredone  & \entryzero          & \entryzero  & \coloredone          & \entryzero  & \entryzero          & \coloredtwo & \entryzero          & \entryzero  & \boldone & \entryzero  & \entryzero          & \entryzero  & \entryzero          & \entryzero  & \entryzero          & \entryzero  & \entryzero          \\
					\matinds[21] & \entryzero          & \entryzero          & \entryzero          & \entryzero          & \entryzero  & \coloredone  & \entryzero  & \entryzero  & \entryzero  & \entryzero  & \entryzero  & \entryzero          & \entryzero  & \entryzero          & \coloredtwo & \entryzero          & \entryzero  & \entryzero          & \entryzero  & \entryzero          & \entryzero  & \entryzero          & \entryzero  & \boldone & \entryzero  & \entryzero          & \entryzero  & \entryzero          & \entryzero  & \entryzero          \\
					\matinds[22] & \entryzero          & \entryzero          & \entryzero          & \entryzero          & \coloredone  & \entryzero  & \entryzero  & \coloredone  & \entryzero  & \entryzero  & \entryzero  & \entryzero          & \entryzero  & \coloredtwo         & \entryzero  & \entryzero          & \coloredtwo & \entryzero          & \entryzero  & \entryzero          & \entryzero  & \entryzero          & \coloredone  & \entryzero          & \entryzero  & \boldone & \entryzero  & \entryzero          & \entryzero  & \entryzero          \\
					\matinds[23] & \entryzero          & \entryzero          & \entryzero          & \coloredone          & \entryzero  & \entryzero  & \coloredtwo & \entryzero  & \entryzero  & \coloredone  & \entryzero  & \entryzero          & \coloredtwo & \entryzero          & \entryzero  & \coloredone          & \entryzero  & \entryzero          & \coloredtwo & \entryzero          & \entryzero  & \coloredone          & \entryzero  & \entryzero          & \coloredtwo & \entryzero          & \entryzero  & \boldone & \entryzero  & \entryzero          \\
					\matinds[24] & \entryzero          & \entryzero          & \coloredone          & \entryzero          & \entryzero  & \entryzero  & \entryzero  & \entryzero  & \entryzero  & \entryzero  & \entryzero  & \entryzero          & \entryzero  & \entryzero          & \entryzero  & \entryzero          & \entryzero  & \entryzero          & \entryzero  & \entryzero          & \entryzero  & \entryzero          & \entryzero  & \entryzero          & \entryzero  & \entryzero          & \entryzero  & \entryzero          & \entryzero  & \boldone \\
					\matinds[25] & \entryzero          & \entryzero          & \entryzero          & \entryzero          & \entryzero  & \entryzero  & \entryzero  & \entryzero  & \entryzero  & \entryzero  & \entryzero  & \entryzero          & \entryzero  & \entryzero          & \entryzero  & \entryzero          & \entryzero  & \entryzero          & \coloredone  & \entryzero          & \entryzero  & \entryzero          & \entryzero  & \entryzero          & \entryzero  & \entryzero          & \entryzero  & \entryzero          & \entryzero  & \entryzero          \\
					\matinds[26] & \entryzero          & \entryzero          & \entryzero          & \entryzero          & \entryzero  & \entryzero  & \entryzero  & \entryzero  & \entryzero  & \entryzero  & \entryzero  & \entryzero          & \entryzero  & \entryzero          & \entryzero  & \entryzero          & \entryzero  & \coloredone          & \entryzero  & \entryzero          & \coloredone  & \entryzero          & \entryzero  & \entryzero          & \entryzero  & \entryzero          & \entryzero  & \entryzero          & \entryzero  & \entryzero          \\
					\matinds[27] & \entryzero          & \entryzero          & \entryzero          & \entryzero          & \entryzero  & \entryzero  & \entryzero  & \entryzero  & \entryzero  & \entryzero  & \entryzero  & \entryzero          & \entryzero  & \entryzero          & \entryzero  & \entryzero          & \coloredone  & \entryzero          & \entryzero  & \coloredtwo         & \entryzero  & \entryzero          & \coloredone  & \entryzero          & \entryzero  & \entryzero          & \entryzero  & \entryzero          & \entryzero  & \entryzero          \\
					\matinds[28] & \entryzero          & \entryzero          & \entryzero          & \entryzero          & \entryzero  & \entryzero  & \entryzero  & \entryzero  & \entryzero  & \entryzero  & \entryzero  & \entryzero          & \entryzero  & \entryzero          & \entryzero  & \coloredone          & \entryzero  & \entryzero          & \entryzero  & \entryzero          & \entryzero  & \entryzero          & \entryzero  & \entryzero          & \coloredone  & \entryzero          & \entryzero  & \entryzero          & \entryzero  & \entryzero          \\
					\matinds[29] & \entryzero          & \entryzero          & \entryzero          & \entryzero          & \entryzero  & \entryzero  & \entryzero  & \entryzero  & \entryzero  & \entryzero  & \entryzero  & \entryzero          & \entryzero  & \entryzero          & \coloredone  & \entryzero          & \entryzero  & \coloredone          & \entryzero  & \entryzero          & \entryzero  & \entryzero          & \entryzero  & \coloredone          & \entryzero  & \entryzero          & \coloredone  & \entryzero          & \entryzero  & \entryzero          \\
					\matinds[30] & \entryzero          & \entryzero          & \entryzero          & \entryzero          & \entryzero  & \entryzero  & \entryzero  & \entryzero  & \entryzero  & \entryzero  & \entryzero  & \entryzero          & \entryzero  & \coloredone          & \entryzero  & \entryzero          & \coloredtwo & \entryzero          & \entryzero  & \coloredone          & \entryzero  & \entryzero          & \coloredone  & \entryzero          & \entryzero  & \coloredtwo         & \entryzero  & \entryzero          & \coloredone  & \entryzero
				\end{bNiceMatrix}_{\textstyle\raisebox{4pt}.}$%
			}
		\end{center}
				
		Then $W(1)_{\z/3\z}$ is easy to get as except the first row, it is a submatrix of $L(1)_{\z/3\z}$.
		Specifically, the matrix obtained by eliminating the first row of $W(1)_{\z/3\z}$ is equal to the matrix obtained by eliminating the 4-th, 14-th, 24-th, 30-th rows and the first, 29-th, 30-th columns from $L(1)_{\z/3\z}$. The matrix $W(1)_{\z/3\z}$ is
		\begin{center}
			\resizebox{\linewidth}{!}{%
				\NiceMatrixOptions{rules/color=[gray]{0.75}, hvlines}
				$\begin{bNiceMatrix}[first-row,first-col]
				\CodeBefore
				\rectangleEntryColor[2-4][2-6]
				\rectangleEntryColor[3-3][3-8]
				\rectangleEntryColor[4-2][4-10]
				\rectangleEntryColor[5-9][5-10]
				\rectangleEntryColor[6-8][6-12]
				\rectangleEntryColor[7-7][7-14]
				\rectangleEntryColor[8-6][8-16]
				\rectangleEntryColor[9-5][9-18]
				\rectangleEntryColor[10-4][10-20]
				\rectangleEntryColor[11-3][11-22]
				\rectangleEntryColor[12-2][12-24]
				\rectangleEntryColor[13-1][13-26]
				\rectangleEntryColor[15-10][15-12]
				\rectangleEntryColor[16-9][16-14]
				\rectangleEntryColor[17-8][17-16]
				\rectangleEntryColor[18-7][18-18]
				\rectangleEntryColor[19-6][19-20]
				\rectangleEntryColor[20-5][20-22]
				\rectangleEntryColor[21-4][21-24]
				\rectangleEntryColor[22-3][22-26]
				\rectangleEntryColor[23-18][23-18]
				\rectangleEntryColor[24-17][24-20]
				\rectangleEntryColor[25-16][25-22]
				\rectangleEntryColor[26-15][26-24]
				\rectangleEntryColor[27-14][27-26]
				\cellcolor[HTML]{CCCCCC}{1-4,2-3,3-2,4-1,5-11,6-13,7-15,8-17,9-19,10-21,11-23,12-25,13-27,14-11,15-13,16-15,17-17,18-19,19-21,20-23,21-25,22-27}
				\Body
					& \matinds[01]  & \matinds[02]  & \matinds[03]  & \matinds[04]  & \matinds[05]  & \matinds[06]  & \matinds[07]  & \matinds[08]  & \matinds[09]  & \matinds[10] & \matinds[11] & \matinds[12] & \matinds[13] & \matinds[14] & \matinds[15] & \matinds[16] & \matinds[17] & \matinds[18] & \matinds[19] & \matinds[20] & \matinds[21] & \matinds[22] & \matinds[23] & \matinds[24] & \matinds[25] & \matinds[26] & \matinds[27]  \\
					\matinds[01]  &   &    &    & \boldone  &    &    &    &    &    &    &    &    &    &    &    &    &    &    &    &    &    &    &    &    &    &    &    \\
					\matinds[02]  &   &    & \boldone  &    &    & \coloredone  &    &    &    &    &    &    &    &    &    &    &    &    &    &    &    &    &    &    &    &    &    \\
					\matinds[03]  &   & \boldone  &    &    & \coloredtwo &    &    & \coloredone  &    &    &    &    &    &    &    &    &    &    &    &    &    &    &    &    &    &    &    \\
					\matinds[04]  & \boldone &    &    &    &    &    &    &    &    & \coloredone  &    &    &    &    &    &    &    &    &    &    &    &    &    &    &    &    &    \\
					\matinds[05]  &   &    &    &    &    &    &    &    & \coloredone  & \coloredtwo & \boldone  &    &    &    &    &    &    &    &    &    &    &    &    &    &    &    &    \\
					\matinds[06]  &   &    &    &    &    &    &    & \coloredone  & \coloredtwo & \coloredone  & \coloredone  & \coloredtwo & \boldone  &    &    &    &    &    &    &    &    &    &    &    &    &    &    \\
					\matinds[07]  &   &    &    &    &    &    & \coloredone  & \coloredtwo & \coloredone  & \coloredtwo & \coloredone  & \coloredtwo & \coloredone  & \coloredtwo & \boldone  &    &    &    &    &    &    &    &    &    &    &    &    \\
					\matinds[08]  &   &    &    &    &    & \coloredone  & \coloredtwo & \coloredone  &    &    &    &    &    &    & \coloredone  & \coloredtwo & \boldone  &    &    &    &    &    &    &    &    &    &    \\
					\matinds[09]  &   &    &    &    & \coloredone  & \coloredtwo & \coloredone  & \coloredone  & \coloredtwo & \coloredone  &    &    &    & \coloredone  & \coloredtwo & \coloredone  & \coloredone  & \coloredtwo & \boldone  &    &    &    &    &    &    &    &    \\
					\matinds[10]  &   &    &    & \coloredone  & \coloredtwo & \coloredone  & \coloredtwo & \coloredone  & \coloredtwo & \coloredone  & \coloredtwo & \coloredone  & \coloredone  & \coloredtwo & \coloredone  & \coloredtwo & \coloredone  & \coloredtwo & \coloredone  & \coloredtwo & \boldone  &    &    &    &    &    &    \\
					\matinds[11] &   &    & \coloredone  & \coloredtwo & \coloredone  &    &    &    &    &    &    & \coloredtwo & \coloredone  & \coloredtwo &    &    &    &    &    &    & \coloredone  & \coloredtwo & \boldone  &    &    &    &    \\
					\matinds[12] &   & \coloredone  & \coloredtwo & \coloredone  & \coloredone  & \coloredtwo & \coloredone  &    &    &    & \coloredtwo & \coloredone  & \coloredtwo & \coloredtwo & \coloredone  & \coloredtwo &    &    &    & \coloredone  & \coloredtwo & \coloredone  & \coloredone  & \coloredtwo & \boldone  &    &    \\
					\matinds[13] & \coloredone & \coloredtwo & \coloredone  & \coloredtwo & \coloredone  & \coloredtwo & \coloredone  & \coloredtwo & \coloredone  & \coloredtwo & \coloredone  & \coloredtwo & \coloredone  & \coloredtwo & \coloredone  & \coloredtwo & \coloredone  & \coloredtwo & \coloredone  & \coloredtwo & \coloredone  & \coloredtwo & \coloredone  & \coloredtwo & \coloredone  & \coloredtwo & \boldone  \\
					\matinds[14] &   &    &    &    &    &    &    &    &    &    & \boldone  &    &    &    &    &    &    &    &    &    &    &    &    &    &    &    &    \\
					\matinds[15] &   &    &    &    &    &    &    &    &    & \coloredone  &    &    & \boldone  &    &    &    &    &    &    &    &    &    &    &    &    &    &    \\
					\matinds[16] &   &    &    &    &    &    &    &    & \coloredone  &    &    & \coloredtwo &    &    & \boldone  &    &    &    &    &    &    &    &    &    &    &    &    \\
					\matinds[17] &   &    &    &    &    &    &    & \coloredone  &    &    &    &    &    &    &    &    & \boldone  &    &    &    &    &    &    &    &    &    &    \\
					\matinds[18] &   &    &    &    &    &    & \coloredone  &    &    & \coloredone  &    &    &    &    &    & \coloredone  &    &    & \boldone  &    &    &    &    &    &    &    &    \\
					\matinds[19] &   &    &    &    &    & \coloredone  &    &    & \coloredtwo &    &    & \coloredone  &    &    & \coloredone  &    &    & \coloredtwo &    &    & \boldone  &    &    &    &    &    &    \\
					\matinds[20] &   &    &    &    & \coloredone  &    &    &    &    &    &    &    &    & \coloredtwo &    &    &    &    &    &    &    &    & \boldone  &    &    &    &    \\
					\matinds[21] &   &    &    & \coloredone  &    &    & \coloredone  &    &    &    &    &    & \coloredtwo &    &    & \coloredtwo &    &    &    &    &    & \coloredone  &    &    & \boldone  &    &    \\
					\matinds[22] &   &    & \coloredone  &    &    & \coloredtwo &    &    & \coloredone  &    &    & \coloredtwo &    &    & \coloredone  &    &    & \coloredtwo &    &    & \coloredone  &    &    & \coloredtwo &    &    & \boldone  \\
					\matinds[23] &   &    &    &    &    &    &    &    &    &    &    &    &    &    &    &    &    & \coloredone  &    &    &    &    &    &    &    &    &    \\
					\matinds[24] &   &    &    &    &    &    &    &    &    &    &    &    &    &    &    &    & \coloredone  &    &    & \coloredone  &    &    &    &    &    &    &    \\
					\matinds[25] &   &    &    &    &    &    &    &    &    &    &    &    &    &    &    & \coloredone  &    &    & \coloredtwo &    &    & \coloredone  &    &    &    &    &    \\
					\matinds[26] &   &    &    &    &    &    &    &    &    &    &    &    &    &    & \coloredone  &    &    &    &    &    &    &    &    & \coloredone  &    &    &    \\
					\matinds[27] &   &    &    &    &    &    &    &    &    &    &    &    &    & \coloredone  &    &    & \coloredone  &    &    &    &    &    & \coloredone  &    &    & \coloredone  &  
				\end{bNiceMatrix}_{\textstyle\raisebox{4pt}.}$%
			}
		\end{center}
		
		We claim $L(1)_{\z/3\z}$ and $W(1)_{\z/3\z}$ are invertible over $\z/3\z$.
		At this point, it may be possible to use Gaussian elimination to show the invertibility as the field $\z/3\z$ is quite simple. However, due to its tediousness, we proceed in another way that one can see it faster.
		For $L(1)_{\z/3\z}$, we want to show for any $1\leq i\leq 30$, there is a column vector of the form \[(\underbrace{0,0,\cdots,0,1}_{i\text{ entries}},*,\cdots,*)\] in the column space whose $i$-th entry is 1 while the entries before it are zero.
		We denote the $i$-th column vector of $L(1)_{\z/3\z}$ as $C_i$.
		For $1\leq i\leq 24$, it is easy to see the entries above the bold faced one can be eliminated by column operations.
		For $25\leq i\leq 30$, the following are the vectors we wanted.
		{\scriptsize
			\begin{alignat*}{10}
				&(0,\cdots,&0&,&1&,&0&,&0&,&0&,&-1&,&1&) &\ =\  &C_3+C_5+C_6+C_{19}+C_{22}-C_{24}-C_{26}-C_{30} \\
				&(0,\cdots,&0&,&0&,&1&,&-1&,&1&,&-1&,&0&) &\ =\  &C_5+C_{21}-C_{23}+C_{25}-C_{27}+C_{28}+C_{29} \\
				&(0,\cdots,&0&,&0&,&0&,&1&,&1&,&-1&,&1&) &\ =\  &C_3+C_6-C_8-C_{15}+C_{17}+C_{24}+C_{25}-C_{26}-C_{27}+C_{28}\\
				&&&&&&&&&&&&&&&&&+C_{29}-C_{30} \\
				&(0,\cdots,&0&,&0&,&0&,&0&,&1&,&-1&,&1&) &\ =\  &-C_1-C_3+C_4-C_6-C_7-C_{17}+C_{20}+C_{22}-C_{23}+C_{24}\\
				&&&&&&&&&&&&&&&&&+C_{25}+C_{27}+C_{28}+C_{30} \\
				&(0,\cdots,&0&,&0&,&0&,&0&,&0&,&1&,&1&) &\ =\  &C_2-C_3-C_4-C_5+C_7-C_8+C_9-C_{10}-C_{11}-C_{13}\\
				&&&&&&&&&&&&&&&&&+C_{14}-C_{17}-C_{18}+C_{21}-C_{22}+C_{23}+C_{26}-C_{27}-C_{29}+C_{30} \\
				&(0,\cdots,&0&,&0&,&0&,&0&,&0&,&0&,&1&) &\ =\  &C_1-C_5-C_6-C_8-C_9-C_{13}-C_{16}+C_{18}+C_{20}-C_{21}\\
				&&&&&&&&&&&&&&&&&-C_{22}+C_{23}+C_{24}+C_{25}+C_{26}+C_{27}-C_{28}\footnotemark
			\end{alignat*}
		}
		\vspace{-11pt}
		\footnotetext{The coefficients of these linear combinations are found by using Macaulay2.}
		
		For $W(1)_{\z/3\z}$, we proceed in a similar way.
		We denote the $i$-th column vector of $W(1)_{\z/3\z}$ as $D_i$.
		For $1\leq i\leq 22$, it is easy to see the entries above the bold faced one can be eliminated by column operations.
		For $23\leq i\leq 27$, the following are the vectors we wanted.
		{\scriptsize
			\begin{alignat*}{9}
				&(0,\cdots,&0&,&1&,&-1&,&1&,&0&,&-1&) &\ =\  &D_3-D_6-D_{16}+D_{18}+D_{19}-D_{20}-D_{21}-D_{25}-D_{26} \\
				&(0,\cdots,&0&,&0&,&1&,&-1&,&1&,&-1&) &\ =\  &-D_5-D_7-D_8-D_{12}-D_{15}+D_{17}+D_{19}-D_{21}+D_{23}-D_{24}+D_{25} \\
				&(0,\cdots,&0&,&0&,&0&,&1&,&1&,&-1&) &\ =\  &D_2+D_5-D_7-D_{14}+D_{16}+D_{23}+D_{24}-D_{25}-D_{26}+D_{27}\\
				&(0,\cdots,&0&,&0&,&0&,&0&,&1&,&-1&) &\ =\  &-D_2+D_3-D_5-D_6-D_{16}+D_{19}+D_{21}-D_{22}+D_{23}+D_{24}\\
				&&&&&&&&&&&&&&&+D_{26}+D_{27} \\
				&(0,\cdots,&0&,&0&,&0&,&0&,&0&,&1&) &\ =\  &D_1-D_2-D_3+D_5+D_6-D_8-D_9-D_{10}+D_{13}+D_{15}\\
				&&&&&&&&&&&&&&&-D_{16}+D_{17}-D_{19}-D_{20}-D_{23}-D_{24}+D_{26}+D_{27}\footnotemark
			\end{alignat*}
		}
		\vspace*{-11pt}
		\footnotetext{The coefficients of these linear combinations are found by using Macaulay2.}
		
		\noindent
		Consequently, $L(1)$ and $W(1)$ are invertible over $\k$, and we get \[\maxdeg X_6\geq 34,\ \deg X_6=33\]
		by \Cref{mainthm_general_m}.
		
		For the second statement, note that $\reg\o_{X_6}\leq \deg X_6-\codim X_6=31$ by \Cref{reg_of_surface}. On the other hand, we have $\reg X_6\geq 34$ by the first statement. This shows $X_6$ should not be 31-normal nor 32-normal.
	\end{proof}
	
	\begin{rem}
		We have proved that there is a homogeneous polynomial $F$ of degree $34$ such that $x_0^{32}x_1^2\in\mon(F)$ in any minimal generating set of $I_{X_6}$. Indeed, using Macaulay2, one can see
		{\footnotesize
			\begin{gather*}
				7\,x_{0}^{32}x_{1}^{2}
				+42\,x_{0}^{29}x_{1}^{2}x_{2}^{3}
				-90\,x_{0}^{26}x_{1}^{2}x_{2}^{6}
				-44\,x_{0}^{23}x_{1}^{2}x_{2}^{9}
				-x_{0}^{20}x_{1}^{2}x_{2}^{12}
				+343\,x_{0}^{27}x_{2}^{2}x_{3}x_{4}^{4}
				\\
				-2058\,x_{0}^{24}x_{2}^{5}x_{3}x_{4}^{4}
				-13230\,x_{0}^{21}x_{2}^{8}x_{3}x_{4}^{4}
				+34398\,x_{0}^{18}x_{2}^{11}x_{3}x_{4}^{4}
				-17913\,x_{0}^{15}x_{2}^{14}x_{3}x_{4}^{4}
				\\
				-1386\,x_{0}^{12}x_{2}^{17}x_{3}x_{4}^{4}
				+2717\,x_{0}^{9}x_{2}^{20}x_{3}x_{4}^{4}
				-555\,x_{0}^{6}x_{2}^{23}x_{3}x_{4}^{4}
				+27\,x_{0}^{3}x_{2}^{26}x_{3}x_{4}^{4}
				+x_{2}^{29}x_{3}x_{4}^{4}
				\\
				+14\,x_{0}^{28}x_{1}x_{4}^{5}
				+280\,x_{0}^{25}x_{1}x_{2}^{3}x_{4}^{5}
				-8118\,x_{0}^{22}x_{1}x_{2}^{6}x_{4}^{5}
				+1599\,x_{0}^{19}x_{1}x_{2}^{9}x_{4}^{5}
				+18695\,x_{0}^{16}x_{1}x_{2}^{12}x_{4}^{5}
				\\
				-10201\,x_{0}^{13}x_{1}x_{2}^{15}x_{4}^{5}
				+175\,x_{0}^{10}x_{1}x_{2}^{18}x_{4}^{5}
				+517\,x_{0}^{7}x_{1}x_{2}^{21}x_{4}^{5}
				-52\,x_{0}^{4}x_{1}x_{2}^{24}x_{4}^{5}
				-x_{0}x_{1}x_{2}^{27}x_{4}^{5}
				\\
				+7\,x_{0}^{24}x_{4}^{10}
				-105\,x_{0}^{21}x_{2}^{3}x_{4}^{10}
				+204\,x_{0}^{18}x_{2}^{6}x_{4}^{10}
				-121\,x_{0}^{15}x_{2}^{9}x_{4}^{10}
				+27\,x_{0}^{12}x_{2}^{12}x_{4}^{10}
				-2\,x_{0}^{9}x_{2}^{15}x_{4}^{10}
			\end{gather*}
		}
		is contained in a minimal generating set of $I_{X_6}$.
	\end{rem}
	
	\begin{rem}
		The linear system we used to construct the map is consisted of four monomials and one binomial in $H^0(\o_{\p^2}(6))$.
		It would be difficult to get a counterexample if the sections were all monomials because the image is a toric variety in that case.
		Note that the projective varieties we constructed are the next simplest projective varieties to those in the sense that all sections in the linear system are monomials except one, which is a binomial.
	\end{rem}
	
	The following proposition shows taking the kernel of a map between polynomial rings over fields commutes with extending the field.
	It shows that the generators of the kernel of a map between the polynomial rings over $\k$ can be obtained by taking the generators of that over $\q$.
	It also enables us to investigate the example over $\q$, instead of $\k$, so that we can use computer algebra systems to compute the projective invariants that are stable under the field extension $\k/\q$.
	
	\begin{prop}\label{primeness_preserved}
		Let $\mathbb{E}/\mathbb{F}$ be a field extension and $I$ be the kernel of the map
		\begin{align*}
			\psi_{\mathbb{F}}:\mathbb{F}[x_0,\cdots,x_r] &\to \mathbb{F}[y_0,\cdots,y_n]\\
			x_i &\mapsto s_i
		\end{align*}
		where all $s_i$'s are homogeneous of the same degree.
		Then the extension $I^e\subseteq\mathbb{E}[x_0,\cdots,x_r]$ of $I$ via the inclusion $\mathbb{F}[x_0,\cdots,x_r]\hookrightarrow\mathbb{E}[x_0,\cdots,x_r]$ is the kernel of the map
		\begin{align*}
			\psi_{\mathbb{E}}:\mathbb{E}[x_0,\cdots,x_r] &\to \mathbb{E}[y_0,\cdots,y_n]\\
			x_i &\mapsto s_i.
		\end{align*}
	\end{prop}
	\begin{proof}
		Let $G$ be a Gr\"{o}bner basis of the ideal \[J=(s_0-x_0,\cdots,s_r-x_r)\subseteq\mathbb{F}[y_0,\cdots,y_n,x_0,\cdots,x_r]\] with respect to the lex order $y_0>\cdots>y_n>x_0>\cdots>x_r$. Then we have $\ker \psi_{\mathbb{F}}=\langle G\cap \mathbb{F}[x_0,\cdots,x_r]\rangle$ (\textit{c.f.} \cite[Theorem 2.4.2]{AL}). The set $G$, considered as a subset of $\mathbb{E}[y_0,\cdots,y_n,x_0,\cdots,x_r]$, is still a Gr\"{o}bner basis of the extension $J^e\subseteq\mathbb{E}[y_0,\cdots,y_n,x_0,\cdots,x_r]$ as the field extension $\mathbb{F}\to\mathbb{E}$ is flat (\textit{c.f.} \cite[Theorem 3.6]{BGS}). Since $G\cap \mathbb{F}[x_0,\cdots,x_r]=G\cap \mathbb{E}[x_0,\cdots,x_r]$, we get $I^e=\ker \psi_{\mathbb{E}}$.
	\end{proof}
	
	In the following example, we investigate some projective invariants, cohomological properties, and geometrical properties of $X_6$ using Macaulay2.
	
	\begin{exmp}\label{example_reg34}
		The surface $X_6$ is the image of the binomial rational map
		\begin{alignat*}{3}
			\varphi:\p^2=\proj\k[y_0,y_1,y_2] &\dashrightarrow &\p^4=\proj \k[x_0,\cdots,x_4]
		\end{alignat*}
		given by $(\o_{\p^2}(6),
		y_1^5y_2,
		y_0^5y_1,
		y_1^6+y_1^3y_2^3,
		y_0^6,
		y_2^5y_0)$.
		We denote these five sections as $s_0,\cdots,s_4$.
		Then its image $X$ is a non-normal surface in $\p^4$ with regularity 34 while the degree is 33. This can be checked over $\q$ by using Macaulay2, and the same must hold over $\k$ by \Cref{primeness_preserved}.
		
		The following is the transposed Betti table of $X_6$.
		
		{\scriptsize
			\begin{center}
			{
				$\begin{array}{c|cccccccccccccccccccc}
					i\backslash j & 0 & \cdots & 7 & 8 & 9 & 10 & 11 & 12 & 13 & 14 & 15 & 16 & 17 & 18 & 19 & 20 & \cdots & 33 \\ \hline
					0             & 1 &        &   &   &   &    &    &    &    &    &    &    &    &    &    &    &        &    \\
					1             &   &        & 1 & 3 &   & 2  & 7  & 16 & 7  & 4  & 2  &    & 3  &    &    & 3  &        & 1  \\
					2             &   &        &   & 1 & 1 & 3  & 14 & 39 & 20 & 11 & 5  & 1  & 8  & 1  &    & 8  &        & 3  \\
					3             &   &        &   &   &   &    & 10 & 23 & 23 & 10 & 4  & 2  & 7  & 2  &    & 7  &        & 3  \\
					4             &   &        &   &   &   &    & 2  & 2  & 9  & 3  & 1  & 1  & 2  & 1  &    & 2  &        & 1
				\end{array}$\\
				Transposed Betti table of $X_6$ with entries $\beta_{i,j}(X_6):=\dim\tor_i(S/I_{X_6},\k)_{i+j}$
			}
			\end{center}
		}
		
		Hence $\reg X_6=\maxdeg I_{X_6}=34> \deg X_6=33$.
		From the cohomological point of view, recall that $\reg \o_{X_6}\leq 31$ and $X_6$ is neither 31-normal nor 32-normal. Indeed, one may check
		\begin{alignat*}{5}
			&h^0(\o_{X_6}(31)) &= P_{X_6}(31) &= 14511,\ \dim (S/I_{X_6})_{31} &= 14509& ,\\
			&h^0(\o_{X_6}(32)) &= P_{X_6}(32) &= 15515,\ \dim (S/I_{X_6})_{32} &= 15514&
		\end{alignat*}
		so that $h^1(\mathcal{I}_{X_6}(31))=2$ and $h^1(\mathcal{I}_{X_6}(32))=1$ in Macaulay2 where $P_{X_6}(n)$ denotes the Hilbert polynomial of $X_6$.
		Partially, the cohomology table of $\o_{X_6}$ with entries $h^i(\o_{X_6}(j-i))$ looks like
		
		{\footnotesize
			\[\begin{array}{c|cccccccccccc}
				i\backslash j & \cdots & -1  & 0   & 1   & 2   & 3   & \cdots & 17   & 18   & \cdots & 34    & \cdots \\ \hline
				2             & \cdots & 255 & 140 & 59  & 14  & 0   & \cdots & 0    & 0    & \cdots & 0     & \cdots \\
				1             & \cdots & 248 & 252 & 260 & 269 & 265 & \cdots & 1    & 0    & \cdots & 0     & \cdots \\
				0             & \cdots & 0   & 1   & 5   & 15  & 35  & \cdots & 3920 & 4462 & \cdots & 17622 & \cdots
			\end{array}\]
		}
		
		\noindent
		so that the regularity of $\o_{X_6}$ is 18.
		
		Note that $X_6$ is a projection of a toric variety $Y_6\subseteq\p^5=\proj\k[x_0,\cdots,x_5]$ of degree $33$ given by the image of the rational map
		\begin{alignat*}{3}
			\p^2          & \dashrightarrow &  & \p^5                                                            \\
			[y_0,y_1,y_2] & \mapsto         &  & [y_1^5y_2,
			y_0^5y_1,
			y_1^6,y_1^3y_2^3,
			y_0^6,
			y_2^5y_0]
		\end{alignat*}
		from the point $q=[0,0,1,-1,0,0]$. One can see $q\notin Y_6$ easily, and $Y_6$ has the following Betti table:
		{\footnotesize
			\begin{gather*}
				\begin{array}{c|ccccc}
					& 0 & 1 & 2  & 3  & 4 \\ \hline
					0  & 1 &   &    &    &   \\
					1  &   &   &    &    &   \\
					2  &   & 2 &    &    &   \\
					3  &   & 1 & 2  &    &   \\
					4  &   &   &    &    &   \\
					5  &   & 1 &    &    &   \\
					6  &   & 1 & 2  &    &   \\
					7  &   & 4 & 12 & 9  & 1 \\
					8  &   & 1 & 3  & 3  & 1 \\
					9  &   & 3 & 10 & 11 & 4 \\
					10 &   &   & 1  & 2  & 1
				\end{array}
			\end{gather*}
		}
		
		To calculate the partial elimination ideals, we take a coordinate change \begin{alignat*}{3}
			\lambda: \p^5          & \to &  & \p^5                                                            \\
			[x_0,x_1,x_2,x_3,x_4,x_5] & \mapsto         &  & [x_0,x_1,x_2+x_3,x_4,x_5,x_2]
		\end{alignat*}
		so that $X_6$ is a projection of $Y_6$ from the point $[0,0,0,0,0,1]$.
		Let $(K_i(I_{Y_6}))_{i\geq 0}$ be the partial elimination ideals of $I_{Y_6}$ with respect to $x_5$.
		As $Y_6\subseteq\p^5$ has regularity 11, it is expected that one of the regularities of the partial elimination ideals of $I_{Y_6}$ is large since \Cref{bound_reg_projection} implies
		\begin{equation}\label{inequ_PEI}
			\reg X_6\leq \max\limits_{1\leq i\leq s-1}(\reg Y_6,\ \reg_S K_i(I_{Y_6})+i+1,\ \reg_S K_s(I_{Y_6})+s-1)
		\end{equation}
		where $s$ is the stabilization number.
		Indeed, using \Cref{PEI_by_GB} and Macaulay2, one can compute \[\reg K_1(I_{Y_6})=32,\ \reg K_2(I_{Y_6})=17,\ \text{and $K_3(I_{Y_6})=(1)$}\]
		so that $s=3$.
		This means the regularity could go bad after the projection due to the first partial elimination ideal $K_1(I_{Y_6})$.
		The locus defined by $K_1(I_{Y_6})$ is quite simple as it is a union of two lines while the ideal itself is complicated in the sense that its regularity is 32.
		Also, this example corresponds to the boundary case of the inequality (\ref{inequ_PEI}).
		It is interesting that the stabilization number is relatively small so that there is no fiber of length $>3$ of the projection.
		
		In fact, the singular loci of $X_6,Y_6$ and the loci of $K_1(I_{Y_6}),K_2(I_{Y_6})$ are
		{\footnotesize
			\begin{alignat*}{6}
				(Y_6)_\sing &= V(x_0,x_1,x_2,x_5)\ &\cup\ &V(x_0,x_1,x_3,x_5) &&\cup V(x_0,x_1,x_3,x_4,x_5-x_2) &\ \subseteq\p^5\\
				(X_6)_\sing &= V(x_0,x_1,x_2) &\cup\ &V(x_0,x_1,x_3) &&\cup V(x_1,x_3,x_4) &\ \subseteq\p^4\\
				V(K_1(I_{Y_6})) &= & &V(x_0,x_1,x_3) &&\cup V(x_1,x_3,x_4) &\ \subseteq\p^4\\
				V(K_2(I_{Y_6})) &= & &V(x_0,x_1,x_2,x_3) &&\cup V(x_1,x_3,x_4) &\ \subseteq\p^4
			\end{alignat*}
		}
		
		\noindent
		where $Y_6$ is the variety coordinate changed by $\lambda$. This also shows $X_6$ is a non-normal surface.
	\end{exmp}
	
	\begin{rem}\label{reg_bound_param_ideal}
		If $X$ is the image of the rational map \[\varphi:\p^n=\proj\k[y_0,\cdots,y_n]\dashrightarrow\p^r=\proj\k[x_0,\cdots,x_r]\]
		given by $(\o_{\p^n}(m),s_0,\cdots,s_r)$, then $\reg X\leq m^{(r+1)2^{n}-1}$ by the recent result of F. Cioffi and A. Conca (\cite[Theorem 1.1]{CC}). In the case of $n=2$, $r=4$, $m=6$ as in \Cref{example_reg34}, we have the bound $\reg X\leq 6^{19}$.
	\end{rem}
	
	In \cite{MP1}, J. McCullough and I. Peeva constructed a threefold counterexample of degree 31 and regularity 38. Recently, J. Choe constructed a threefold counterexample of degree 18 and regularity 18 in \cite{C2}.
	Here, we give a surface counterexample of degree 11 and regularity 11 which is the counterexample with the smallest regularity among those we are aware of.
	
	\begin{exmp}\label{example_reg11}
		Let \[\varphi:\p^2=\proj\k[y_0,y_1,y_2]\dashrightarrow\p^4=\proj\k[x_0,\cdots,x_4]\] be a rational map given by $(\o_{\p^2}(4),y_0^3y_1,y_0^2y_1^2,y_1^2y_2^2,y_0y_2^3+y_1^3y_2,y_0^2y_1y_2+y_1^4)$.
		Using Macaulay2, one can see its image is a non-normal surface in $\p^4$ with the Betti table
		{\footnotesize
			\[\begin{array}{c|ccccc}
				& 0 & 1  & 2  & 3  & 4 \\ \hline
				0      & 1 &    &    &    &   \\
				\vdots &   &    &    &    &   \\
				4      &   & 1  &    &    &   \\
				5      &   & 13 & 21 & 8  &   \\
				6      &   & 4  & 12 & 12 & 4 \\
				7      &   & 2  & 5  & 4  & 1 \\
				8      &   &    &    &    &   \\
				9      &   &    &    &    &   \\
				10     &   & 1  & 3  & 3  & 1
			\end{array}\]
		}
		
		\noindent
		so that it has regularity 11. The degree is also 11 in this case.
	\end{exmp}
	
	Similarly, there is a threefold counterexample $X\subseteq\p^5$ with $\reg X=\deg X=11$ that is not a cone.
	
	\begin{exmp}\label{example_reg11_threefold}
		Let \[
		\varphi:\p^3=\proj\k[y_0,\cdots,y_3] \dashrightarrow \p^5=\proj\k[x_0,\cdots,x_5]
		\] be a rational map given by $(\o_{\p^3}(3),y_0y_1^2,y_0y_2y_3,y_1y_2^2,y_2^3,y_0^2y_1+y_2^2y_3,y_0^2y_2+y_1^3)$.
		Using Macaulay2, one can see its image is a non-normal threefold $X\subseteq\p^5$ with the Betti table
		{\footnotesize
			\[\begin{array}{c|cccccc}
				& 0 & 1  & 2  & 3  & 4 & 5 \\ \hline
				0      & 1 &    &    &    &   &  \\
				\vdots &   &    &    &    &   &  \\
				4      &   & 1  &    &    &   &  \\
				5      &   & 8 & 11 & 3  &   &  \\
				6      &   & 18  & 53 & 56 & 25 & 4 \\
				7      &   & 3  & 9  & 10  & 5 & 1 \\
				8      &   &    &    &    &   &  \\
				9      &   &    &    &    &   &  \\
				10     &   & 1  & 3  & 3  & 1 &
			\end{array}\]
		}
		
		\noindent
		so that it has regularity 11. The degree is also 11. As the projective dimension of $S_X$ is 5, the depth of $S_X$ is 1. Consequently, $X$ is not a cone.
	\end{exmp}
	
	\section{Questions}
	Although we were able to find some counterexamples to the regularity conjecture, we could not find one that is normal. Note that P. Mantero, L. E. Miller, and J. McCullough proved that the threefold suggested by J. McCullough and I. Peeva in \cite[Example 4.7]{MP1} is regular in codimension one in \cite[Theorem A]{MMM} while it is not normal. To the authors' knowledge, normal counterexamples are not known currently.
	Note that being normal is a weaker condition than both of being smooth and being projectively normal.
	Hence it would be interesting if one can answer the following question:
	\begin{question}
		Is there a counterexample to the Eisenbud-Goto conjecture that is normal?
	\end{question}
	
	It is possible to construct a threefold that is regular in codimension one also using a rational map between projective spaces: the image of the rational map
	\[
	\varphi:\p^3=\proj\k[y_0,\cdots,y_3] \dashrightarrow \p^5=\proj\k[x_0,\cdots,x_5]
	\]
	given by $(\o_{\p^3}(3),y_0^3,y_0^2y_1,y_0y_2^2,y_1^2y_3,y_0y_2y_3+y_1^3,y_1^2y_2+y_2y_3^2)$ for example.
	In this case, the regularity is 19, the degree is 15, and the singular locus is a union of five lines and a space curve. However, for surfaces, we could not find any counterexample that is regular in codimension one. Hence for the case of surfaces, before considering normal surfaces, we raise the following question:
	
	\begin{question}
		Is there a counterexample to the Eisenbud-Goto conjecture that is a surface with isolated singularities?
	\end{question}
	
	Note that the regularity conjecture holds for normal surfaces with mild singularities due to W. Niu (\cite{N1}).
	
	\section*{Acknowledgement}
	We are grateful to A. Conca, J. McCullough, and I. Peeva for their careful reading and helpful comments.
	A. Conca pointed out there is a bound on the regularity of prime ideals with polynomial parametrisations as in \Cref{reg_bound_param_ideal}.
	J. McCullough informed us about the study on the singular loci of counterexamples to the regularity conjecture (\cite{MMM}). We are also grateful to the anonymous referee for detailed review with valuable suggestions.
	
	\section*{Statements and Declarations}
	There are no competing interests to declare.


\begin{thebibliography}{XXXX}
		\bibitem[AL]{AL} W. Adams and P. Loustaunau, An Introduction to Grobner Bases, American Mathematical Society, 1994.
		
		\bibitem[BGS]{BGS} D. Bayer, A. Galligo, and M. Stillman, Gr\"{o}bner bases and extension of scalars, Computational Algebraic
		Geometry and Commutative Algebra (Cortona, 1991), Cambridge Univ. Press, 1993, 198–215.
		
		\bibitem[C1]{C1} G. Castelnuovo, Sui multipli di una serie lineare di gruppi di punti appartenente ad una curva algebrica, Rend. Circ. Mat. Palermo 7 (1893), 89-110.
		
		\bibitem[C2]{C2} J. Choe, Castelnuovo-Mumford regularity of unprojections and the Eisenbud-Goto regularity conjecture, preprint, arXiv:2206.06151, 2022.
		
		\bibitem[CC]{CC} F. Cioffi and A. Conca, Regularity of primes associated with polynomial parametrisations, Ann. Sc. Norm. Super. Pisa Cl. Sci. XXIV (2023), 1149-1154.
		
		\bibitem[CCG]{CCG} L. Chiantini, N. Chiarli, and S. Greco, Bounding Castelnuovo–Mumford regularity for varieties with good general projections, J. Pure Appl. Algebra 152 (2000), 57-64.
		
		\bibitem[CS]{CS} A. Conca and J. Sidman, Generic initial ideals of points and curves, J. Symbolic Comput. 40 (2005), 1023–1038.
		
		\bibitem[EG]{EG} D. Eisenbud and S. Goto, Linear free resolutions and
		minimal multiplicity, J. Algebra 88 (1984), 89–133.
		
		\bibitem[G]{G} M. Green, Generic initial ideals, Six Lectures on Commutative Algebra, Birkh\"{a}user, 1998, 119-186.
		
		\bibitem[GLP]{GLP} L. Gruson, R. Lazarsfeld, and C. Peskine, On a theorem	of Castelnuovo, and the equations defining space curves, Invent. Math. 72 (1983), 491–506.
		
		\bibitem[K1]{K1} A. G. Kouchnirenko, Polyèdres de Newton et nombres de Milnor, Invent. Math. 32 (1976), 1-32.
		
		\bibitem[K2]{K2} S. Kwak, Castelnuovo regularity for smooth subvarieties of dimensions 3 and 4, J. Algebraic
		Geom. 7 (1998), 195–206.
		
		\bibitem[K3]{K3} S. Kwak, Castelnuovo–Mumford regularity bound for smooth threefolds in $\p^5$ and extremal
		examples, J. Reine Angew. Math. 509 (1999), 21–34.
		
		\bibitem[KNV]{KNV} S. Kwak, H. Nguyen, and T. Vu, Algebraic invariants of projections of varieties and partial elimination ideals, J. Algebra 586 (2021), 973-1013.
		
		\bibitem[KP]{KP} S. Kwak and J. Park, A bound for Castelnuovo-Mumford regularity by double point divisors, Adv. Math. 364 (2020), 107008.
		
		\bibitem[L]{L} R. Lazarsfeld, A sharp Castelnuovo bound for smooth surfaces, Duke Math. J. 55 (1987), 423-429.
		
		\bibitem[M2]{M2} D. R. Grayson and M. E. Stillman, Macaulay 2, a software system for research in algebraic geometry,
		\url{http://www.math.uiuc.edu/Macaulay2/}.
		
		\bibitem[MMM]{MMM} P. Mantero, L. E. Miller, and J. McCullough, Singularities of Rees-like algebras, Math. Z. 297 (2021), 535-555.
		
		\bibitem[MP1]{MP1} J. McCullough and I. Peeva, Counterexamples to the Eisenbud-Goto regularity conjecture, J. Amer. Math. Soc. 31 (2018), 473-496.
		
		\bibitem[MP2]{MP2} J. McCullough and I. Peeva, The Regularity Conjecture for prime ideals in polynomial rings, EMS Surv. Math. Sci. 7 (2020), 173–206.
		
		\bibitem[N1]{N1} W. Niu, Castelnuovo–Mumford regularity bounds for singular surfaces, Math. Z. 280 (2015), 609-620.
		
		\bibitem[N2]{N2} A. Noma, Generic inner projections of projective varieties and an application to the positivity of double point divisors, Trans. Amer. Math. Soc. 366 (2014), 4603-4623.
		
		\bibitem[NP1]{NP1} W. Niu and J. Park, A Castelnuovo–Mumford regularity bound for
		scrolls, J. Algebra 488 (2017), 388-402.
		
		\bibitem[NP2]{NP2} W. Niu and J. Park, A Castelnuovo-Mumford regularity bound for threefolds with rational singularities, Adv. Math. 401 (2022), 108320.
		
		\bibitem[PS]{PS} I. Peeva and B. Sturmfels, Syzygies of codimension 2 lattice ideals, Math. Z. 229 (1998), 163-194.
		
		\bibitem[R]{R} Z. Ran, Local differential geometry and generic projections of threefolds, J. Differential Geom. 32 (1990), 131-137.
		
		\bibitem[T]{T} S. Telen, Introduction to toric geometry, preprint, arXiv:2203.01690, 2022.
	\end{thebibliography}
\end{document}